\documentclass[11pt,notitlepage,twoside,a4paper]{amsart}
 \usepackage{amsfonts}
\usepackage[latin1]{inputenc}
\usepackage[cyr]{aeguill}
\usepackage{xspace}
\usepackage[polutonikogreek,english]{babel}

\usepackage{amsmath,amssymb,enumerate}

\usepackage{epsfig,fancyhdr,color}

\usepackage{amssymb}
\usepackage{amsmath,amsthm}
\usepackage{latexsym}
\usepackage{amscd}
\usepackage{psfrag}
\usepackage{graphicx}
\usepackage[latin1]{inputenc}
\usepackage[all]{xy}
\usepackage[mathcal]{eucal}

\definecolor{NoteColor}{rgb}{1,0,0}

\renewcommand{\textsc}{\textcolor{red}}

\newtheorem*{theorem 1}{\rm\bf Proposition 1}
\newtheorem*{theorem 2}{\rm\bf Proposition 2}

\theoremstyle{definition}

\theoremstyle{remark}

\def\interieur#1{\mathord{\mathop{\kern 0pt #1}\limits^\circ}}

\title[W. P. Thurston and French mathematics]{W. P. Thurston and French mathematics}
\author[Fran\c cois Laudenbach, Athanase Papadopoulos et al.]{Fran\c cois Laudenbach and Athanase Papadopoulos,
 \\
 {}
 \\
 {\smaller\rm with contributions by William Abikoff,
  Norbert A'Campo, 
  \\
  Pierre Arnoux,  Michel Boileau, 
  Albert Fathi, 
David Fried, 
\\
Gilbert Levitt, Valentin Po\'enaru, 
Harold Rosenberg,  
\\
Francis Sergeraert, Vlad Sergiescu and Dennis Sullivan}}

\begin{document}

\maketitle
 
\date{\today}


 \begin{abstract}
 We give a general overview of the influence of William Thurston on the French mathematical school and we show how some of the major problems he solved are rooted in the French mathematical tradition. At the same time, we survey some of Thurston's major results and their impact. The final version of this paper will appear in the Surveys of the European Mathematical Society.

  \medskip
    \medskip
    \noindent  AMS classification: 01A70; 01A60; 01A61; 57M50;  57R17; 57R30
      
  \medskip
  \medskip
   \noindent   Keywords: William P. Thurston: Geometric structures; Hyperbolic structures; Haefliger structures;  low-dimensional topology; foliations; contact structures; history of French mathematics.
                       
               \end{abstract}          
                        
         \part{}
                \section{Prologue} 
 Seven years have passed since Bill Thurston left us, but his presence is felt every day in the minds of a whole community of mathematicians who were shaped by his ideas and his completely original way of thinking about mathematics.

In 2015-2016 a two-part celebration of Thurston and his work was published in the {\em Notices of the AMS}, edited by Dave Gabai and Steve Kerckhoff with contributions by several of Thurston's students and other mathematicians who were close to him \cite{Gabai}.  Among the latter was our former colleague and friend Tan Lei, who passed away a few years later, also from cancer, at the age of 53.   One of Tan Lei's last professional activities was the thesis defense of her student J\'er\^ome Tomasini in Angers, which took place on December~5, 2014  and for which  
Dylan Thurston served on the committee.  Tan Lei was suffering greatly at that time; her disease was diagnosed a few days later. She confessed to Fran\c cois Laudenbach that when she received the galley proofs of the article on Thurston which included her contribution, she was sad to realize that there was hardly any mention of Thurston's influence on the French school of mathematics, and in particular on the Orsay group.  In several subsequent phone calls with Laudenbach she insisted that this story needed to be told.  Our desire to fulfill her wish was the motivation for the present article.  We asked for the help of a number of colleagues who were involved in the mathematical activity ``in the tradition of Thurston" that took place (and continues to take place) in France.  Most of these colleagues responded positively and sent us a contribution, either in the form of a short article or as an email, everyone in his personal style and drawing from his own memory.  Sullivan, after we showed him our article, proposed to include a text of his which previously appeared in the \emph{Notices}, considering that his text is complementary to ours and would naturally fit here. We happily accepted. All these contributions are collected in the second part of this article. We have included them as originally written, although there are a few minor discrepancies in some details pertaining to the description or the precise order of the events.   

In the first part of this article we have tried to give a general overview of the influence of Thurston on the French mathematical community.  At the same time we show how some of the major problems he solved have their roots in the French mathematical tradition.

Saying that Thurston's ideas radically changed the fields of low-dimensional geometry and topology and had a permanent effect on related fields such as geometric group theory and dynamics, is stating the obvious. On a personal level we remember vividly his generosity, his humility and his integrity.  

We view this article as another general tribute to Thurston, besides being a chronicle of the connection of his work with that of French mathematicians.  It is not a survey of Thurston's work---such a project would need several volumes, but it is a survey on the relation between Thurston and French mathematics. At the same time, we  hope that it will give the reader who is unfamiliar with Thurston's \oe uvre an idea of its breadth.

\section{\emph{Vita}}

It seemed natural to us to start with a short \emph{Vita}.
William Paul  Thurston was born in Washington DC on Oct. 30, 1946 and he died on August 21, 2012.
His father was an engineer and his mother a housewife. He went to Kindergarten in Holland and then entered the American educational system. 
 
Thurston's  parents were wise enough to let him choose the college he wanted to attend after he completed his secondary education.  In 1964, he entered New College, a small private college in Sarasota (Florida). The college was newly founded. 
In a commencement address he gave there in 1987 (that was twenty years after his graduation), Thurston recalls that he landed in New College because, by chance, he read an advertisement for that college,  amidst the large amount of  literature which used to fill the family mail box every day. In that advertisement, two statements of educational principles attracted his attention. The first one was: ``In the final analysis, each person is responsible for her or his education," and the second one: ``The best education is the collaboration/conflict between two first-class minds." For him, he says, these statements looked like ``a declaration of independence, of a sort of freedom from all the stupidities of all the schools that [he has]  sort of rebelled against all [his] life."
Thurston recalled that in grade school he was spending his time daydreaming and (maybe  with a small amount of exaggeration) his grades were always Cs and Ds. In Junior high school, he became more rebellious with his teachers and was traumatised by the fact that these teachers, because they were in a position of authority, were supposed never to make  mistakes.
He kept fighting, as a child and then as a teenager, with the American educational system and he could not adapt to it  until he went to college.

In New College, there was a focus on independent study, and writing a senior thesis was one of the important requirements. Thurston was interested in the foundations of mathematics and this motivated him to write a thesis on intuitionist topology. The title of the thesis, submitted in 1967, is:  ``A Constructive Foundation for Topology."
He was particularly attracted by intuitionism, not only as a topic for a senior thesis, but he thought he might become an intuitionist logician. This is why,   when he entered the University of California, Berkeley, he approached Alfred Tarski, the charismatic logician and mathematics teacher there, and asked him to be his advisor. Tarski told the young Thurston that Berkeley was not a good place for intuitionism. Thus, Thurston went instead to topology, the other theme of his senior thesis.

Thurston, especially in his last years, used to actively participate in discussions on mathematics and other science blogs, and reading his contributions gives us some hints on his way of thinking. The key word that is recurrent in his writings about his view on the goal of mathematics is ``understanding".  In a post dated May 17, 2010, he writes: 
\begin{quote}\small
 In my high school yearbook, I put as my goal ``to understand", and I've thought of that as summing up what drives me.
My attitude toward the demarcation problem originated I think from childhood games my siblings and I used to play, where one of us would say something obviously implausible about the world, as if psychotic, and the others would try to trip up the fantasy and establish that it couldn't be right.  We discovered how difficult it is to establish reality, and I started to think of these battles as futile. People of good will whose thinking is not confused and muddled or trapped in a rut can reach a common understanding. In the absence of good will or clarity, they do not, and an external criterion or external referee does not help. 

[\ldots]  I also used to think I would switch to biology when I reached the age of 35 or 40, because I was very drawn to the challenge of trying to understand life. It didn't happen.
 \end{quote}

Thurston received his PhD in 1972; his thesis is entitled {\it ``Foliations of Three-Manifolds which are Circle Bundles''}. He was appointed full professor at Princeton in 1974  (at the age of 27), 
and he remained there for almost 20 years.  He then returned to California where he joined the University of  California, Berkeley (1992--1997), acting also as the director of MSRI, then the University of California, Davis (1996--2003). After that, he moved to Cornell where he spent the last 9 years of his life.

  In the mid-1970s, Thurston formulated a conjecture on the geometry of 3-manifolds which is the analogue in that dimension of the fact that any surfaces (2-manifold) carries a metric of constant curvature. (The 2-dimensional statement is considered as a form of the uniformization theorem.) The conjecture, called Thurston's conjecture,  says tat any 3-dimensional manifold  can be decomposed in a canonical way into pieces such that each piece carries one of eight types of geometric structures that became known as Thurston geometries. The 3-dimensional geometrization conjecture is wider than the Poincar\'e conjecture (actually a question formulated by Poincar\'e in 1904 at the very end of the {\it Cinqui\`eme compl\'ement \`a l'Analysis Situs}), saying that a simply-connected 3-dimensional manifold without boundary is homeomorphic to the 3-sphere.\footnote{Poincaré made the following comment: 
  ``Cette question nous entra\^{i}nerait trop loin'', that is, {\it This question would take us too far}.}  Thurston's geometrization conjecture was proved in 2002 by Grigory Perelman. At a symposium held in Paris in 2010 celebrating the proof of the Poincar\'e conjecture, Thurston recalls: ``At a symposium on Poincar\'e in 1980, I felt emboldened to say that the geometrization conjecture put the Poincar\'e conjecture into a fuller and more constructive context." He then adds: ``I expressed confidence that the geometrization conjecture is true, and I predicted it would be proven, but whether in one year or 100 years I could not say---I hoped it would be within my lifetime. I tried hard to prove it. I am truly gratified to see my hope finally become reality."

 In 1976 Thurston was awarded the Oswald Veblen Geometry Prize of the American Mathematical Society for his work on foliations, in 1982 the Fields Medal and in 2012 the Leroy P. Steele Prize of the American Mathematical Society.

Thurston introduced a new way of communicating and writing mathematics. He had a personal  and unconventional idea on what mathematics is about and why we do mathematics, and he tried to share it. On several occasions,  he insisted that mathematics does not consist of definitions, theorems and proofs, but of ways of seeing forms and patterns, of internalizing and imagining the world, and of thinking and understanding certain kinds of phenomena. He was attached to the notion of mathematical community. After he finished college, he realized the existence of such a community, and this  appeared to him like a revelation. In 2012, in his response to the Leroy P. Steele Prize, he declared: ``I felt very lucky when I discovered
the mathematical community---local, national and international---starting
in graduate school." In a post on mathoverflow (October 30, 2010), he wrote:
\begin{quote}\small
Mathematics only exists in a living community of mathematicians that spreads understanding and breathes life into ideas both old and new. The real satisfaction from mathematics is in learning from others and sharing with others. All of us have clear understanding of a few things and murky concepts of many more. There is no way to run out of ideas in need of clarification. The question of who is the first person to ever set foot on some square meter of land is really secondary. Revolutionary change does matter, but revolutions are few, and they are not self-sustaining---they depend very heavily on the community of mathematicians. 
\end{quote}
Thurston saw that a school of thought sharing his geometric vision was gradually growing. 
In the same response to the Leroy P. Steele Prize, he declared: 
\begin{quote}\small 
I used to feel that there was certain knowledge
and certain ways of thinking that were unique to me. It is very satisfying to
have arrived at a stage where this is no longer true---lots of people have picked
up on my ways of thought, and many people have proven theorems that I once
tried and failed to prove.''
\end{quote}

\section{Foliations}

  Thurston made his first visit to the Orsay department of mathematics in 1972, the year he obtained his PhD. He was invited by Harold Rosenberg, who was a professor there. Orsay is a small city situated south of Paris, about 40 minutes drive from Porte d'Orl\'eans (which is one of the main south entrances to Paris) with the usual traffic jam.
   The Orsay department of mathematics was very young;  it was created in 1965 as part of a project to decentralize the University of Paris, and in 1971 it had  become  part of the newly founded Universit\'e de Paris-Sud.  Thurston was working on foliations and France was, at that time, the world center for this topic.  Specifically, the theory was born some two-and-a-half decades earlier in Strasbourg, where a strong group of topologists had formed around Charles Ehresmann, including Ren\'e Thom, Georges Reeb, Andr\'e Haefliger and Jean-Louis Koszul.  Reeb's PhD thesis, entitled
\emph{Propri\'et\'es topologiques des vari\'et\'es feuillet\'ees} \cite{Reeb}
which he defended in 1948, may safely be considered
to be the birth certificate of foliation theory. Reeb described there the first example of a foliation of the 3-sphere, answering positively a question asked by Heinz Hopf in 1935. This question was communicated to Reeb by his mentor, Ehresmann (Reeb mentioned this several times). On the other hand, examples of the use of foliations of surfaces can be traced back to the early works on cartography by Ptolemy and others before him who searched for mappings of the 2-sphere onto a plane where the foliation of the sphere by parallels or by longitudes is sent to foliations of a planar surface satisfying certain \emph{a priori} conditions (circles, ellipses, straight lines, lines that interpolate between circles near the North Pole and straight lines near the South Pole, etc.)

    At the beginning of the 1970s the theory of foliations was a hot research topic among topologists and dynamicists at Orsay.  The results that Thurston obtained during his graduate studies and the  years immediately following (ca. 1970--1975) constituted a striking and unforeseeable breakthrough in the field.  In this period of 4 or 5 years, he solved all the major open problems on foliations, a development which eventually led to the disappearance of the Orsay foliation group.
  
 During his 1972 visit, Thurston lectured on his version of the $h$-principle---as it is now called---for foliations of 
codimension greater than 1, and in particular on his result saying that an arbitrary field of 2-planes on a manifold of dimension at least four is homotopic to a smooth integrable one, that is, a field tangent to a foliation. Several young researchers interested in foliations attended the lectures, including Robert Roussarie, Robert Moussu, Norbert A'Campo, Michel Herman and Francis Sergeraert.

 The Godbillon-Vey invariant (GV)  had been born the year before in Strasbourg,\footnote{Claude Godbillon,  a former student of Reeb, was a professor there, and  the discovery was made during a visit to him by  
Jacques Vey (1943--1979) who was a young post-doc.}  
but it was even unknown whether the invariant could be nonzero.   
Informed of this during a meeting at Oberwolfach which took place on May 23-29, 1971, Roussarie immediately found an example of a foliation 
with a non-zero GV invariant, 
namely the horocyclic foliation on a compact quotient of $SL_2(\mathbb R)$.  Shortly after, Thurston proved the much stronger result saying that there exists a family of foliations whose GV take all possible real values. He also gave a geometric interpretation of the GV class of a foliation $F$ 
 as a ``helical wobble of the leaves of $F$\,".  The paper was published in 1972 \cite{Thurston1972}.
 This was one of the first papers that Thurston published on foliations (in fact, it was his first paper on the subject after his thesis, the latter of which remains unpublished). Thurston sent the preprint of his paper to Milnor. On November 22, 1971, Milnor responded with a 5 page letter. We have reproduced here the first page of that correspondence.
  Thurston spent the next academic year at Princeton's Institute for Advanced Study, at the invitation of Milnor. It is interesting to note that the hyperbolic plane already appears in Thurston's paper as a central object. Orbifolds also appear in the background (Thurston calls them ``surfaces having a number of isolated corners, with metrics of constant negative curvature everywhere else").  Shortly thereafter  Sullivan, who was working at IH\'ES in Bures-sur-Yvette (a few minutes walk from Orsay), gave another interpretation of the GV invariant using a notion of linking number for currents; this appeared
 in his 1976 paper \emph{Cycles for the dynamical study of foliated manifolds and complex manifolds} \cite{Sullivan-Cycles}.

\begin{figure}[htbp]
\centering
\includegraphics[scale=.8]{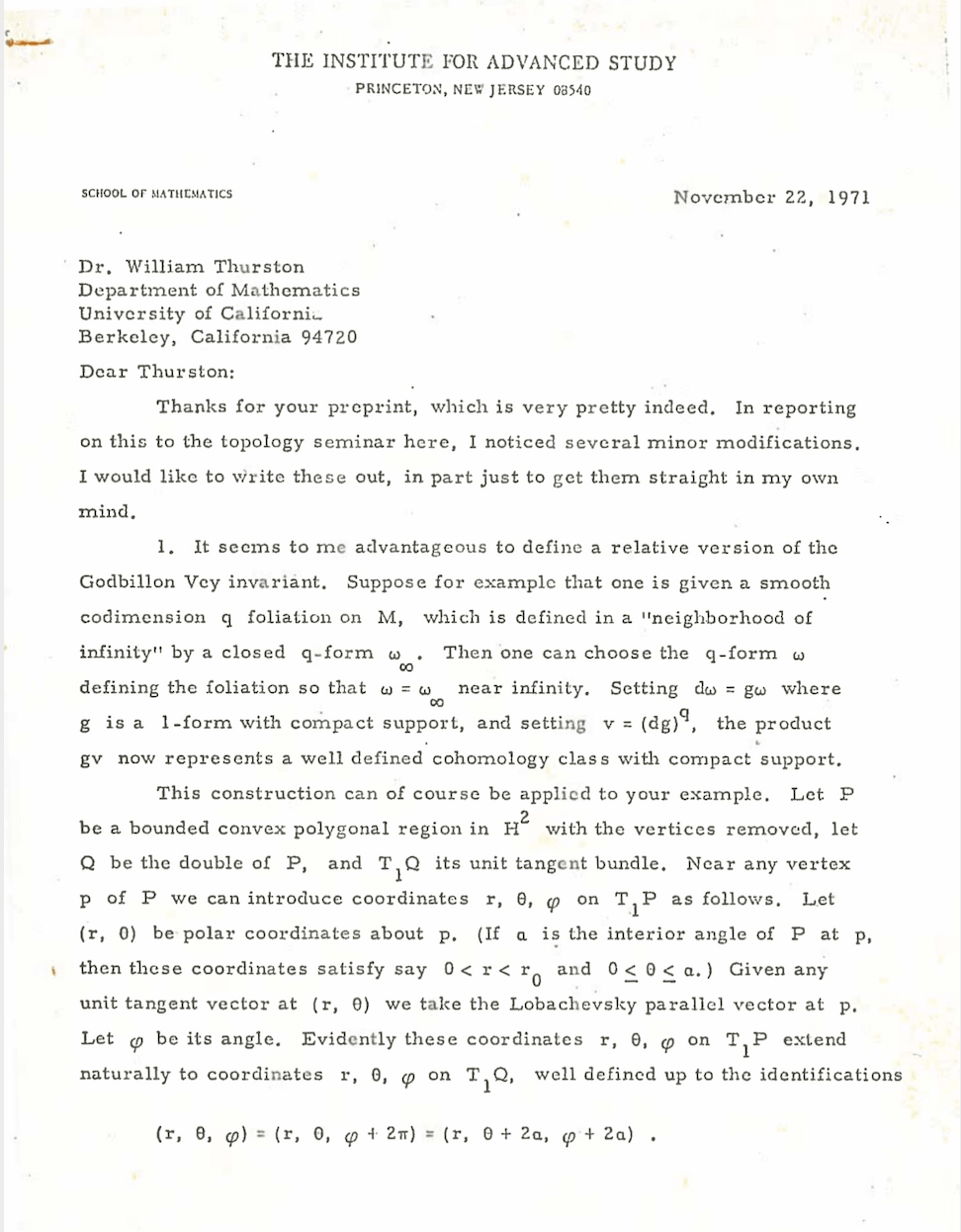}
\caption{The first page of a letter from Milnor to Thurston. Courtesy of H. Rosenberg.} \label{Milnor}
\end{figure}

  During his 1972 stay in France, Thurston also visited  Dijon, and he then went to Switzerland to participate in a conference on foliations at Plans-sur-Bex, a village in the Alps.\footnote{In an obituary article on Michel Kervaire, written by S. Eliahou, P. de la Harpe,
 J.-C. Hausmann and M. Weber and
published in the \emph{Gazette des Math\'ematiciens} \cite{Eliahou},  the authors write that Haefliger knew that a young student from Berkeley had obtained remarkable results on foliations. The news that a  week on foliations would be organized spread rapidly, and the organisers were almost refusing people because the housing possibilities were limited. There were not enough chairs in the classroom of the village school for all the participants but they managed to find one for him. It goes without saying that this student was William Thurston.} A'Campo recalls that at that time, Thurston was already thinking about hyperbolic geometry in dimension two. He asked Thurston how he came to know about this subject, and Thurston's answer was that his father first told him about it.

 Thurston's most remarkable result during the period that followed his Orsay visit is probably the proof of the existence of a $C^\infty$ codimension-one foliation on any closed manifold with zero Euler characteristic. The result was published in 1976 \cite{T-existence}, and it solved one of the main conjectures in the field. Before that, there were only particular examples of foliations of spheres and some other particular manifolds. Sullivan recalls that the first new foliation of spheres after Reeb's was constructed by Lawson in Bahia in 1971, and Verjovsky helped in that. Other particular examples of foliations of special manifolds, sometimes restricted to a single dimension, were constructed by A'Campo,
Durfee, Novikov and Tamura.  On the other hand  Haefliger had proved a beautiful and influential result saying that there is no codimension-one real-analytic foliation on a sphere of any dimension  
 \cite{Haefliger1956}.  This emphasizes the fact
that unlike manifolds, a smooth one-dimensional non-Hausdorff manifold---as is, in general, the space of leaves of a foliation---does not carry any analytic structure. It is interesting to note here that each of Haefliger and Reeb, at a conference in Strasbourg in 1944--1955, presented a paper on non-Hausdorff manifolds considered as quotient spaces of foliations \cite{Haefliger1954, Reeb1954}.

  In his paper \cite{Thurston-stability} published in 1974, Thurston obtained a generalization of the so-called Reeb stability theorem. This theorem, proved by Reeb twenty-eight years earlier in his thesis \cite{Reeb}, says that if a codimension-one foliation of a compact manifold has a two-sided compact leaf with finite fundamental group, then all the leaves of the foliation are  diffeomorphic.  Thurston showed that in the case of a $C^1$ foliation one can replace the hypothesis that the compact leaf has finite fundamental group by the much weaker one saying that the first real cohomology group  of the leaf is zero. He also gave a counterexample in the case where the smoothness condition is not satisfied. As a corollary, he showed that there are many manifolds with boundary that do not admit foliations tangent to the boundary.

  At the same time, Thurston proved in \cite{thurston-cmh} a series of breakthrough results on codimension-$k$
Haefliger structures when $k>1$. Such a structure, introduced (without the name) by Haefliger in his thesis \cite{Haefliger1958}, is a generalization of a foliation: it is an $\mathbb R^k$-bundle over an $n$-dimensional manifold equipped with a codimenson-$k$ foliation transverse
to the fibers of the bundle. (The normal bundle to a foliation is naturally equipped with such a structure.)

Two years later, Thurston solved the tour de force case of codimension-one: every hyperplane field is homotopic to the tangent plane field of a $C^\infty$-foliation \cite{T-existence}. In the same paper, he writes  that ``the theory of analytic foliations still has many unanswered questions."

Thurston gave a talk at the 1974 ICM (Vancouver) whose title is \emph{On the construction and classification of foliations}.
The proceedings of this congress  
contain a short paper (3 pages) \cite{Thurston-Vancouver}   in which he states his major results. The definition he gives of a foliation is unusual, but it  delivers the meaning of the object defined better than any formal definition: ``A foliation is a manifold made out of a striped fabric---with infinitely thin stripes, having no space between them. The complete stripes, or `leaves', of the foliation are submanifolds; if the leaves have codimension $k$, the foliation is called a 
codimension-$k$ foliation".

After Thurston proved his series of results on foliations, the field stagnated.  In his article \emph{On proof and progress in mathematics} \cite{BAMS1994} Thurston wrote: 
  \begin{quote} \small
 Within a couple of years, a dramatic evacuation
of the field started to take place. I heard from a number of mathematicians that they were giving or receiving advice not to go into foliations---they were saying
that Thurston was cleaning it out. People told me (not as a complaint, but as a compliment) that I was killing the field. Graduate students stopped studying foliations.
  \end{quote}
  
 Thurston was never proud of this, and it was not his intention to kill the subject.  In the same article he explained that on the contrary he was sorry about the fact that many people abandoned the field and he thought the new situation arose out of a misunderstanding, by a whole group of mathematicians, of the state of the art of foliations. As a matter of fact, in the introduction to his  
 1976 paper \cite{T-existence} which supposedly killed the field, Thurston declares that further work on the subject is called for, especially using the geometrical methods of his predecessors  and unlike his own method which, according to his words, is local and has the disadvantage that it is hard to picture the foliations constructed.

Thurston might be pleased to know that some problems on foliations that he was interested in have been revived recently.  For instance, G. Meigniez  recently revisited Thurston's ideas and obtained new unexpected results \cite{meigniezJDG}. In particular, he showed  that there exist minimal, $C^\infty$, codimension-one foliations on every closed connected manifold of dimension at least 4 whose Euler characteristic is zero. Since by definition every leaf of a minimal foliation is dense, this proves that there is no generalization to higher dimensions of Novikov's 3-dimensional compact leaf theorem \cite{Novikov} (1965).

\section{Contact geometry}

At the same time he was working on foliations, Thurston obtained a number of important results on contact geometry.  This subject is close to the theory of foliations even though, by definition, a contact structure is very different from a foliation: it is a  hyperplane field  that is maximally far from being integrable, that is, of being tangent to a foliation. 

When Thurston started working on contact structures, the theory, like that of foliations, was already well developed in France. As in the case of the theory  of foliations, a group of topologists in Strasbourg had been working on contact structures since the second half of the 1940s, under the guidance of Ehresmann.
In particular Wu Wen Tsu, who was a young researcher in Strasbourg at the time, published two papers in 1948 in which he studied the existence of contact structures (as well as almost-complex structures) on spheres and sphere bundles. His motivation came from some problems on characteristic classes, a subject on which  he was then competing with Thom. 

In the 1970s, contact geometry was one of the favorite objects of study at the mathematics institute in Strasbourg. Robert Lutz, Jean Martinet and others were working on it under the guidance of Reeb. In his thesis defended in 1971, Lutz proved that every homotopy class of co-orientable plane fields 
on the 3-sphere contains a contact structure. In the same year, Lutz and Martinet, improving techniques used by Lutz in his thesis, showed that every closed orientable 3-manifold supports a contact structure. In 1975, Thurston came in. He published a paper with H. E. Winkelnkemper \cite{TW} 
which gives an amazingly short proof (less than one page) of the result of Lutz and Martinet by using the so-called {\it open-book decomposition} theorem of Alexander. Several years later, and still in dimension three, Emmanuel Giroux obtained a much more difficult result, namely, 
any contact structure is {\it carried} by an open-book decomposition.

Together with Eliashberg, Thurston later developed the notion of \emph{confoliation} in dimension three and  techniques of approximating smooth foliations by contact structures, thus further strengthening the links between the two subjects.

 In higher dimensions, contact structures are much more complicated (hypotheses are needed for existence) and the complete picture is still not well understood. Nevertheless, Giroux over the course of several years obtained 
a generalization to higher dimensions of his theorem that we mentioned above in dimension three, using ideas originating in Donaldson's
asymptotically holomorphic sequences of sections adapted to contact structures
by Ibort, Martinez-Torres and Presas.

Thurston continued to think about contact structures. In contrast with the theory of foliations,  the theory of contact structures is still extremely active in France.

      A further close relation between  contact structures, Thurston and French mathematics is given by  the  Bennequin--Thurston\footnote{The story of this discovery (which is not joint work, but work in parallel, by Thurston and Bennequin) is intricate, and we were not able to reconstruct the exact chronology. Thus, we decided to follow the alphabetical order.} invariant of a Legendrian knot, which describes its amount of coiling.  This in turn gave rise to the Bennequin--Thurston number of a knot, which maximizes  the Bennequin--Thurston invariant over all Legendrian representatives; these invariants were found independently by Thurston and by Daniel Bennequin.        An inequality conjectured by Thurston and proved by Bennequin in his thesis (1982)  has since become known  as the Bennequin--Thurston inequality.
      
 Thurston had a very personal way of explaining contact structures (and the same can be said regarding almost any topic that he talked about). The chapter titled ``Geometric manifolds" of his monograph \emph{Three-dimensional geometry and topology} \cite{Thurston1} contains a section dedicated to contact structures. This comes between the section on bundles and connections and the one on the eight geometries. Thurston spends several pages trying to give an intuitive picture of contact structures, because, he says, ``they give an interesting example of a widely occurring pattern for manifolds that is hard to see until your mind and eyes have been attuned."  On p. 172 of this monograph, he writes:  
 \begin{quote}\small
 You can get a good physical sense for the contact structure on the tangent circle bundle of a surface by thinking about ice skating, or bicycling. A skate that is not scarping sideways describes a Legendrian curve in the tangent circle bundle to the ice. It can turn arbitrarily, but any change of position is in the direction it points. Likewise, as you cycle along, the direction of the bicycle defines a ray tangent to the earth at the point of contact of the rear wheel. Assuming you are not skidding, the rear wheel moves in the direction of this ray, and this motion describes a Legendrian curve in the tangent bundle of the earth.
 
 Young children are sometimes given bicycles with training wheels, some distance off to the side of the rear wheel. The training wheel also traces out a Legendrian curve---in fact, for any real number $t$, the diffeomorphism $\phi_t$ of $\mathbb{R}^2\times S^1$ that takes a tangent ray a signed distance $t$ to the left of itself is a contact automorphism. The training wheel path is the image of the rear wheel path under such a transformation. Note that this transformation applied to curves in the plane often creates or removes cusps. The same thing happens when you mow a lawn, if you start by making a big circuit around the edge of the lawn and move inward [...]
 \end{quote}
We have quoted this long passage because it is characteristic of Thurston's style, that of providing mental images with analogies borrowed from the real world. At the same time, it may be appropriate to recall Thurston's warning  that ``one person's clear mental image is another person's intimidation" \cite{BAMS1994}.

 \section{Hyperbolic geometry, surfaces and 3-manifolds}
        
Among the important classical subjects that Thurston revived one finds the study of hyperbolic structures on surfaces and 3-manifolds.
This brings us to the second sensational piece of  Thurston's work  that had a long-lasting impact on the Orsay group of geometry,  namely, his work on surface mapping class groups, Teichm\"uller spaces, and the geometry and topology of 3-manifolds.

\begin{figure}[htbp]
\centering
\includegraphics[scale=.66]{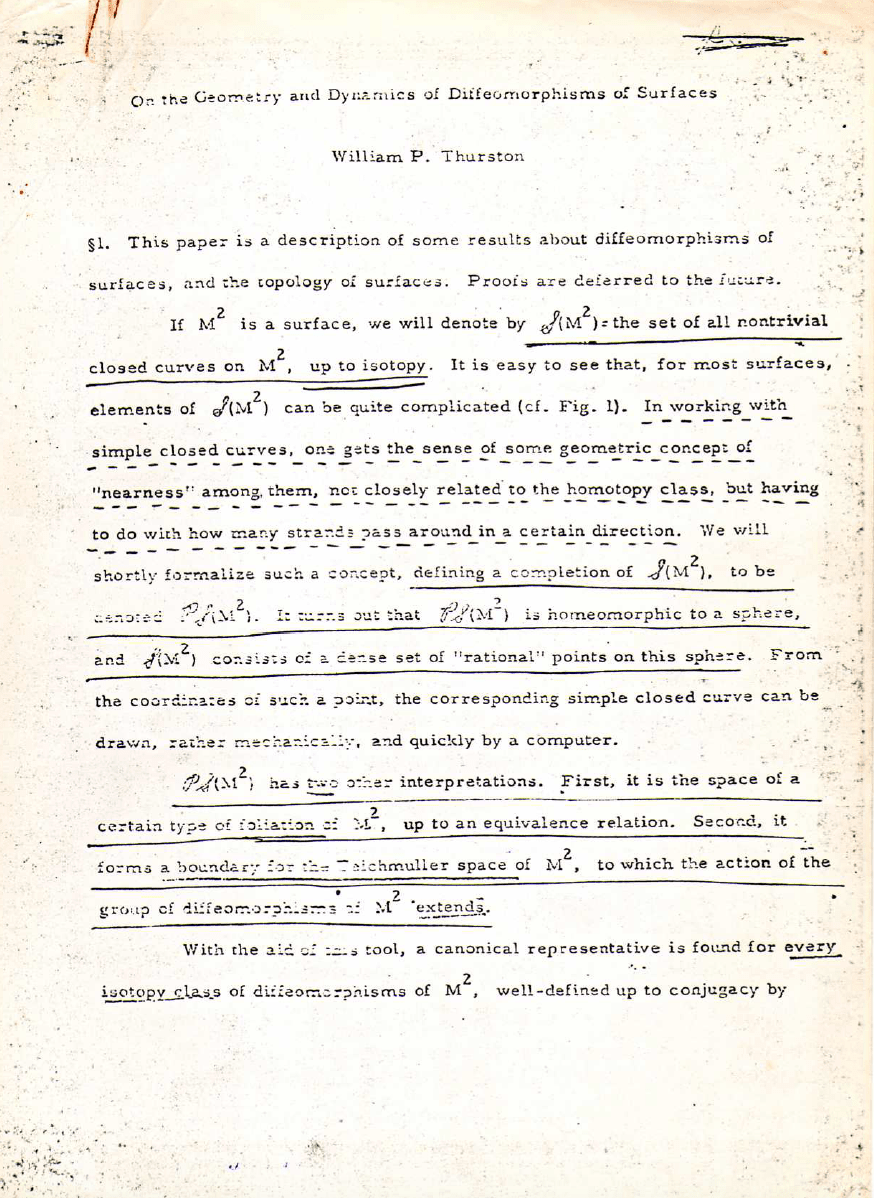}
\caption{The first page of the first document that was the basis of the book ``Travaux de Thurston sur les surface".} \label{first}
\end{figure}

\begin{figure}[htbp]
\centering
\includegraphics[scale=.26]{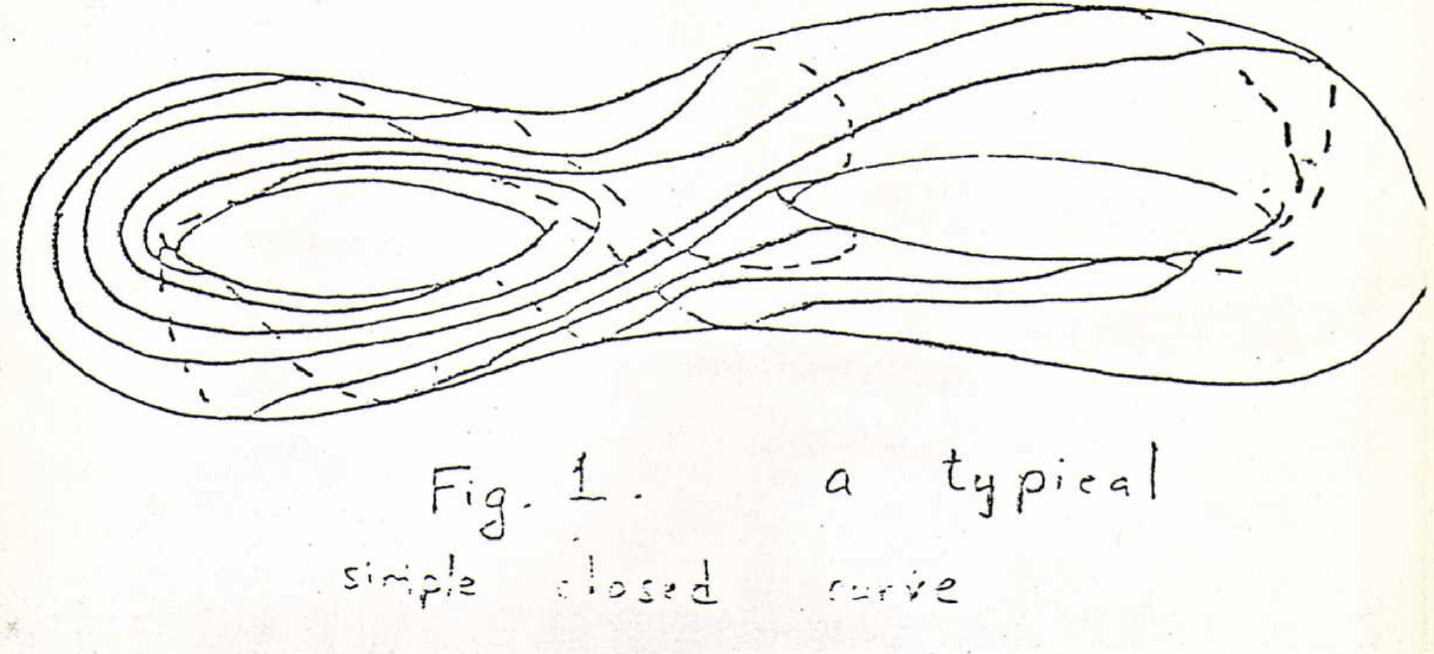}
\quad\includegraphics[scale=.23]{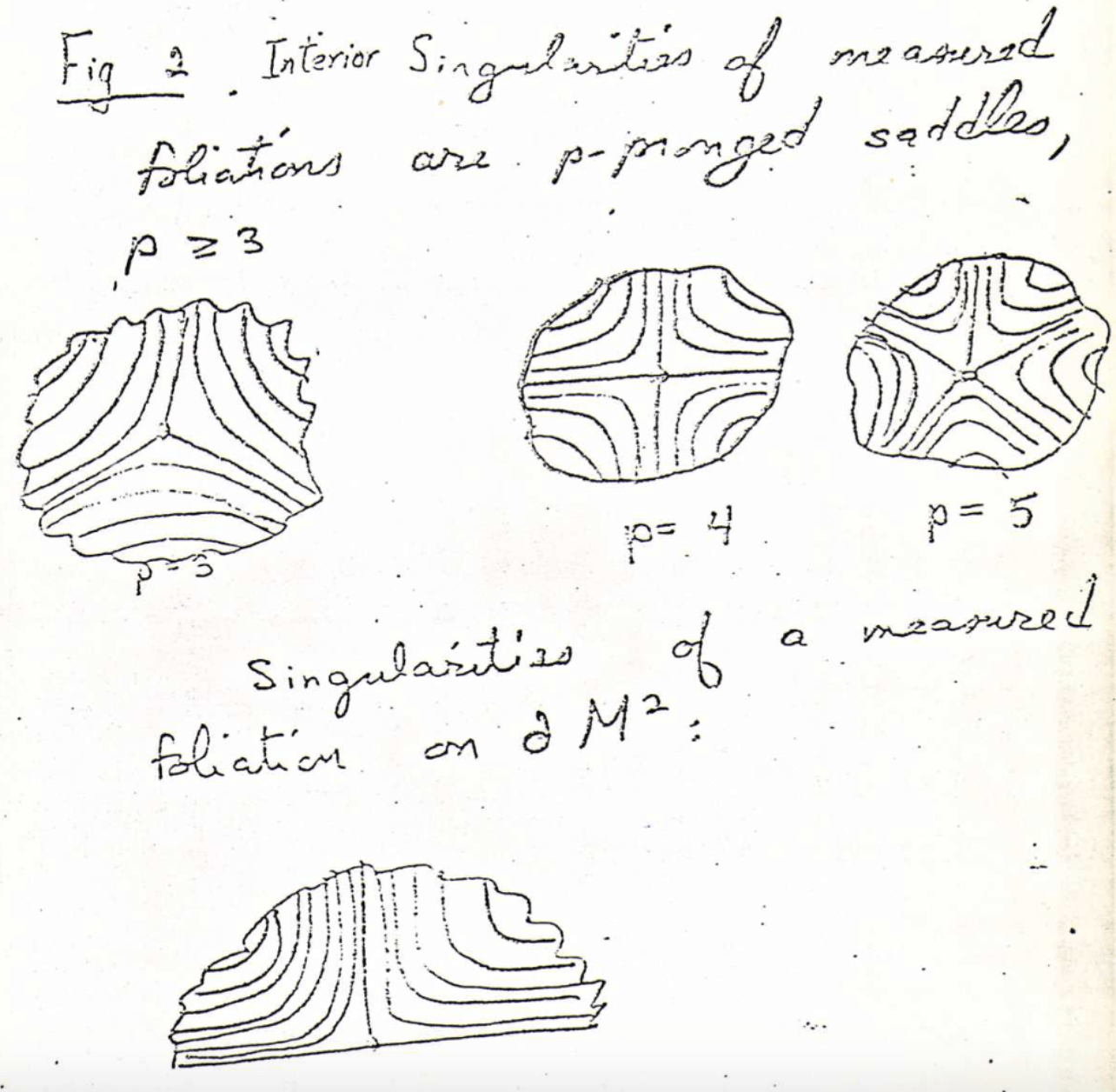}\quad
\includegraphics[scale=.30]{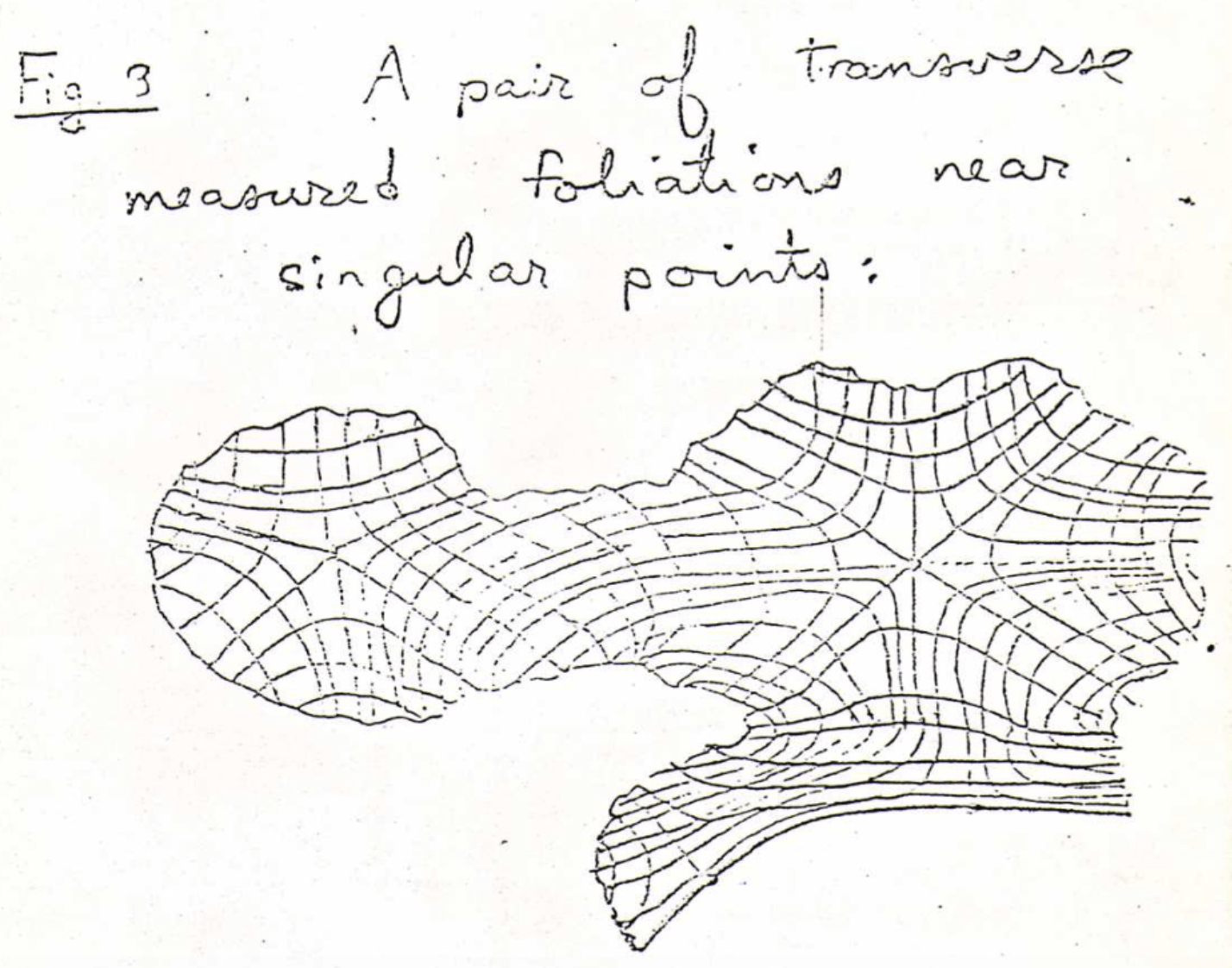}
\caption{A long simple curve; singular foliations; a pair of transverse foliations.}\label{drawings}
\end{figure}

The idea to launch the Orsay seminar known as \emph{Travaux de Thurston sur les surfaces} came from Valentin Po\'enaru.
 At the 1976 ``Autumn cocktail" of IH\'ES, he showed up with a set of notes by Thurston, whose first page is shown in Figure \ref{first}. The notes contained the outline (definitions, pictures, and the statements of results) of what became known later as Thurston's theory of surfaces.\footnote{In reality, Thurston had been thinking about simple closed curves on surfaces ever since he was a student in Berkeley; see the Second Story in the recollections by Sullivan in Part II of the present paper.} These notes were published several years later in the Bulletin of the AMS under the title \emph{On the geometry and dynamics of diffeomorphisms of surfaces} \cite{Thurston-FLP}. Figure \ref{drawings} shows the drawings by Thurston at the end of his preprint.

 The Orsay seminar on the works of Thurston on surfaces took place in the academic year 1976-77. 
Thurston never attended this seminar.  At some point  the seminar members were stuck with the problem of
  gluing the space of projective measured
foliations to Teichm\"uller space as a boundary, and they asked him for assistance.  At that time there was no email, only postal mail.  They finally managed to work out a complete
answer only when Thurston attended another session of the seminar at Plans-sur-Bex in 1978.\footnote{In the obituary article on Kervaire  cited previously \cite{Eliahou}, the authors list the foreign participants at that meeting as A. Connes,
D. Epstein, M. Herman, D. McDuff, J. Milnor, V. Po\'enaru, L. Siebenmann, D. Sullivan and
W. Thurston.}

The final result  of this seminar was published in book form in 1979 in \emph{Ast\'erisque}, and inspired many young topologists and 
geometers. One feature of the book was a first course in hyperbolic geometry, a topic for which, at that time, there were very few reference books. Some important notions that were introduced by Thurston shortly after the seminar took place were missing,
in particular {\it train tracks} and {\it geodesic laminations}. We learned about them after Thurston's manuscript on surfaces reached Orsay, in the mimeographed notes of his Princeton 1978 course on the topology and geometry of 3-manifolds. These were certainly new and fundamental concepts, but the theory of surfaces worked quite well without them.

  Thurston's notes on 3-manifolds arrived in France in installments and their importance was immediately realized by those who had been following his work. They were xeroxed and bound chapter by chapter, in several dozen copies. They were made available in the secretarial office of the Orsay topology research group, to all the members of the group, but also to the other mathematicians who were curious about the theory. In those days the Orsay topology seminar attracted a large number of mathematicians from all over France,  and for many of them, the secretarial office was a necessary passage; they had to sign papers there. In this way  the whole community of French topologists became aware of these notes.
 
 During the year that the \emph{Orsay seminar} took place and during  several years after that, a significant number of foreign mathematicians visited Orsay and gave  lectures and courses on topics related to Thurston's ideas. These included Bill Abikoff,  Lipman Bers, David Epstein, David Fried, John Morgan, Peter Shalen, Mike Shub, and many others. Joan Birman came with two students, John McCarthy and J\'ozef Przytycki. Bob Penner also visited as a student.  Most of these mathematicians continued to maintain strong relations with their French colleagues. A'Campo and Po\'enaru gave graduate courses that were attended by many students and colleagues. Starting in the early 1980s, several doctoral dissertations were defended at Orsay and Paris 7 in which the authors had benefitted 
 from the book issued from the Orsay seminar, the mimeographed notes and the courses given at Orsay.  Among the early graduates were Gilbert Levitt (Th\`ese d'\'Etat 
1983),  Claude Danthony (PhD thesis 1986) and Athanase Papadopoulos (Th\`ese d'\'Etat 1989), who worked on surfaces. Francis Bonahon (Th\`ese d'\'Etat 1985), Michel Boileau (Th\`ese d'\'Etat 1986) and Jean-Pierre Otal (Th\`ese d'\'Etat 1989) worked on 3-manifolds.  Otal later wrote a book which became a standard reference for Thurston's hyperbolization theorem for fibered 3-manifolds.

 After obtaining their doctoral degrees, these young geometers obtained jobs in various places in France---it was  a period of ``decentralization" for mathematics appointments in France, especially at CNRS. One seminar was organized in Strasbourg  by Morin and Papadopoulos, under the name GT3 (in honor of Thurston's notes on the Geometry and Topology of 3-manifolds).  Wolpert, Floyd, Mosher, Epstein, Bowditch, Fried, Gabai and others visited this seminar. McCarthy, Oertel and Penner were long-term visitors in Strasbourg. All of them were involved in Thurston type geometry, and Thurston's results on surface geometry and 3-manifolds were discussed extensively.  The seminar still runs today.  A series of results on Thurston's asymmetric metric on Teichm\"uller space,  after those of Thurston's foundational  paper  \emph{Minimal Stretch maps between hyperbolic
surfaces} \cite{Thurston1985} (1985), were obtained by Strasbourg researchers.  They concern the boundary behavior of stretch lines, the action of the mapping class group on this metric, the introduction of an analogous metric on Teichm\"uller spaces of surfaces with boundary (the so-called arc metric), and there are several other results. Thurston's paper \cite{Thurston1985} was (and is still) considered as being difficult to read, although it uses only material from classical geometry and first principles.  This shows ---if proof is needed--- that profound and difficult mathematics remains the one that is based on simple ideas. Other groups of topologists influenced by Thurston were formed in Marseille (Lustig, Short), Toulouse (Boileau, Otal), and various other places in France. 
   Bonahon moved to the US after he proved a major result that was conjectured by Thurston in his Princeton lecture notes, namely, that the ends of a hyperbolic 3-manifold whose fundamental group is isomorphic to that of a closed surface are geometrically tame \cite{Bonahon}. (The notion of tameness was introduced by Thurston.)  A'Campo went to teach in Basel and introduced several young mathematicians to Thurston-type geometry and topology. His student Walter Br\"agger gave a new proof of Thurston's version of Andreev's theorem, see \cite{Braegger}.

      \section{Holomorphic dynamics}

Dennis Sullivan was a major promoter of Thurston's ideas in France, and he was probably the person who best understood their  originality and implications. For  more than twenty years Sullivan ran a seminar at IH\'ES on topology and dynamics. Recurrent themes at that seminar were Kleinian groups (discrete isometry groups of hyperbolic 3-space), a subject whose foundations were essentially set by Poincar\'e, and holomorphic dynamics, another subject rooted in French mathematics, namely in the works of Fatou and Julia, revived 60 years later, by Adrien  Douady and John Hubbard, in the early 1980s, but preceded sometime in the late 1970s by Milnor and Thurston who developed their so-called kneading theory for the family of maps $x\mapsto x^2+c$. The topics discussed in Sullivan's seminar also included the geometry of 3-manifolds, deformations of Kleinian groups and their limit sets, pleated surfaces, positive eigenfunctions of the Laplacian, quasiconformal mappings, and one-dimensional dynamics. Thurston's ideas were at the forefront, and Sullivan spent years explaining them.

In 1982, while Sullivan was running his seminar on holomorphic dynamics, Douady gave a course on the same subject   at Orsay. At the same time, Sullivan established his dictionary between the iteration theory of rational maps and the dynamics of Kleinian groups. 
  
 In the same year, Sullivan was the first to learn from Thurston about his theorem characterizing  postcritically finite rational maps of the sphere, that is, rational maps whose forward orbits of critical points are eventually periodic. The proof of this theorem, like the proofs of several of Thurston's big theorems, uses  a fixed point argument for an action on a Teichm\"uller space.  Specifically, Thurston associated to a  self-mapping of the sphere which is
 postcritically finite a self-map of the Teichm\"uller space  of the sphere with some points deleted (the postcritical set). The rational map in the theorem is then obtained through an iterative process as a fixed point of the map on Teichm\"uller space. 
 
In addition to the map on Teichm\"uller space, the proof of Thurston's theorem involves hyperbolic geometry, the action of the mapping class group of the punctured sphere on essential closed curves and the notions of  invariant laminations.  All of these notions form the basis for a beautiful analogy between the  ideas and techniques used in the proof of this theorem and those used in the proof of Thurston's classification of mapping classes  of surfaces, and this correspondence is an illustration of the fact that mathematics, for Thurston, was a single unified field.

      Thurston circulated several versions of a manuscript in which he gave all the ingredients of the proof of his theorem, but the manuscript was never finished.\footnote{The first version, circulated in 1983, carries the title \emph{The combinatorics of iterated rational maps}; the subsequent versions known to the authors of the present article do not carry any title.} A proof of this theorem following Thurston's outline was written by Douady and Hubbard. A first version was circulated in  preprint form in 1984 and the paper was eventually published in 1993 \cite{DH}.

Three years later, Sullivan published a paper in which he gave the proof of a long-standing question formulated by Fatou and Julia \cite{Sullivan1985} (1985). The result became known as the No-wandering-domain Theorem. It says that every component of the Fatou set of a rational function is eventually periodic. A fundamental tool that was introduced by Sullivan in his proof is that of quasiconformal mappings, one of the main concepts in classical Teichm\"uller theory. These mappings became a powerful tool in the theory of iteration of rational maps. It is interesting that Sullivan, in his paper \cite{Sullivan1985}, starts by noting that the perturbation of the analytic dynamical system $z\mapsto z^2$ to $z\mapsto z^2+az$ for small $a$ strongly reminds one of Poincaré's perturbations of Fuchsian groups $\Gamma\subset \mathrm{PSL}(2,\mathbb{R})$ into quasi-Fuchsian groups in $\mathrm{PSL}(2,\mathbb{C})$ where the Poincaré limit set changes from a round circle to a non-differentiable Jordan curve, and that Fatou and Julia, the two founders of the theory of iteration of analytic mappings, were well aware of the analogy with Poincaré's work. He then writes: ``We continue this analogy by injecting the modern theory of quasiconformal mappings into the dynamical theory of iteration of complex analytical mappings."

             Thurston's  theorem, together with Sullivan's dictionary, now constitute the two most fundamental results in the theory of iterations of rational maps.

      In his PhD thesis, defended under Thurston in Princeton in 1985, Silvio Levy obtained several applications of Thurston's theorem, including a condition for the existence of a {\em mating} of two degree-two polynomials that are postcritically finite \cite{Levy}. The notion of mating of two polynomials of the same degree was introduced in 1982 by Douady and Hubbard. The idea was to search for a rational self-map of the sphere that combines the dynamical behavior of the two polynomials.   Levy, in his  thesis, formulated the question of mating in a more combinatorial way, and using Thurston's characterization of rational maps  was able to give a necessary and sufficient condition for the existence of a mating of two postcritically finite degree-two polynomials in terms of their associated laminations. This result solved a question formulated in several precise forms by Douady in his Bourbaki seminar \cite[Questions 11 and 12]{D-B}.
            At the same time, also in his thesis, Levy established connections between Thurston's geometric approach and Douady-Hubbard's more analytical approach to the subject of iteration of rational maps. In particular, he established the relation between Thurston's invariant laminations and the so-called Hubbard trees that were introduced by Douady and Hubbard in the context of degree-two polynomials. Both notions arise from identifications that arise on the boundary of the unit disc when it is sent by a Riemann mapping to the complement (in the Riemann sphere) of the so-called filled Julia set of a polynomial, in the case where this set is connected.

       Tan Lei's thesis, which she defended in 1986 at Orsay under the supervision of  Douady, is in some sense the French counterpart of  Levy's thesis. It uses  Thurston's theorem in an essential way, but instead of laminations Tan Lei works with Hubbard trees. A criterion that Douady and Hubbard formulated in \cite[III. 3]{D-B} gives a necessary condition for the existence of a rational function realizing the mating  two degree-two polynomials. Tan Lei provided sufficient conditions for this to happen, at the same time giving  a more precise form of the criterion found by  Levy.        
      
In her paper \emph{Branched coverings and cubic Newton maps} \cite{Tan} Tan Lei applies Thurston's theory of postcritically finite branched coverings of the sphere to a new family of maps. 
Specifically, she studies the dynamics of a class of degree-3 rational maps that arise in Newton's method for approximating the roots of a cubic polynomial. She introduces the notion of a postcritically finite cubic Newton map and  investigates the question of whether branched coverings of the sphere  are equivalent (in the sense of Thurston) to such a map.  The problem of understanding and giving precise information about the roots of a complex polynomial is one of those basic mathematical questions which Thurston was always interested in.  

In 2011, about a year and a half before his death, Thurston posted a thread  on math overflow concerning the intersection of the convex hull of level sets  $\{z \vert Q(z) = w\}$  for a polynomial $Q$. He writes: ``By chance, I've discussed this question a bit with Tan Lei; she made some nice movies of how the convex hulls of level sets vary with $w$. (Also, it's fun to look at their diagrams interactively manipulated in Mathematica). If I get my thoughts organized I'll post an answer." Thurston never had a chance to post the answer.

 Motivated by this question, Tan Lei wrote an article 
 with Arnaud Ch\'eritat  in the French electronic journal \emph{Images des math\'ematiques} dedicated to the popularization of mathematics. In this article,  they first present a classical result known as the Gauss--Lucas Theorem, saying that the convex hull of the roots of any polynomial $P$ of degree  at least one contains the roots of its derivative $P'$.\footnote{Gauss implicitly used this result in 1836, while he formulated the problem of locating the zeros of the derivative of a polynomial in a mechanical way:  he showed that these zeros (provided they are distinct from the multiple roots of the polynomial), are the equilibrium positions of the field of force generated by identical particles placed at the roots of the polynomial itself and where each particle generates a force of attraction which satisfies the inverse proportional distance law \cite{Gauss}. Lucas, in 1874, published a mechanical proof of the same theorem, while he was unaware of Gauss's work \cite{Lucas}. (Gauss's notes, published later in his \emph{Collected Works}, were still poorly known.)} Note that the roots of the derivative are the critical points of the original polynomial.
Tan Lei and Ch\'eritat present a result of Thurston which gives a complete geometric picture of the situation:
\emph{Let $P$ be a non-constant polynomial. Let $F$ be a half-plane bounded by a support line of the convex hull of the roots of the derivative $P'$ of $P$ and not containing this convex hull and let $c$ be a root of $P'$ contained on this support line. Then there is a connected region contained in $F$ on which $P$ is bijective and whose interior is sent by $P$ onto a plane with a slit along a ray starting at $P(c)$.} Tan Lei and Ch\'eritat  gave the details of Thurston's proof that avoids computations, and they provided the computer movies that Thurston talked about in his post. Their article constitutes a tribute to Bill Thurston; it was published less than 3 months after his death. A more detailed version, including two more authors, Yan Gao and Yafei Ou, was later published in the Comptes Rendus \cite{Cheritat}. 

 Over the years, Douady's courses on holomorphic dynamics at Orsay were attended by a number of students and also by more senior mathematicians, including John Hubbard, Pierette Sentenac, Marguerite Flexor, Tan Lei, Pierre Lavaurs,  Jean Ecalle, S\'ebastien Godillon, Arnaud Ch\'eritat, Ricardo P\'erez-Marco, Xavier Buff, and Jean-Christophe Yoccoz.\footnote{It is probably under Douady's influence that the topology research unit called \emph{\'Equipe de topologie} at the University of Orsay was replaced by a unit called \emph{\'Equipe de topologie et dynamique} which until today forms one of the five research units at the mathematics department there.} 
      Tan Lei, in her tribute to Thurston in \cite{Gabai}, writes that he never stopped thinking about iterations of rational maps. She gives a lively description of her conversations and email exchanges with him on this subject in 2011 and 2012, the last two years of his short life.

\begin{figure}[htbp]
\centering
\includegraphics[scale=.70]{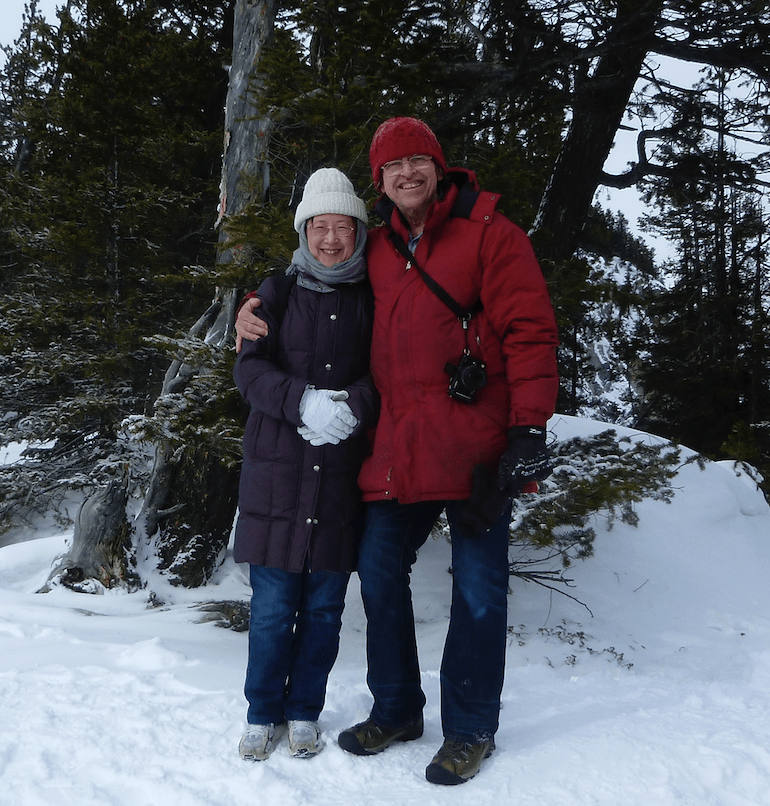}
\caption{Bill Thurston and Tan Lei, Banff, February 2011. Photo courtesy of A. Ch\'eritat and H. Rugh }\label{Thurston-Tan}
\end{figure}

In the realm of conformal geometry, Thurston  introduced the subject of discrete conformal mappings, and in particular the  idea of discrete Riemann mappings.
 In 1987, Sullivan, together with Burton Rodin, proved an important conjecture of Thurston on approximating a Riemann mapping using circle packings
 \cite{SR}. Colin de Verdi\`ere, motivated by Thurston's work, proved the first variational principle for circle packings \cite{CV}.

              \section{Corrugations}

   After foliations and contact structures, let us say a word about corrugations.
     
 From the very beginning of his research activity, Thurston was interested in immersion theory. 
 This was shortly after the birth of the {\it h-principle} in Gromov's 1969 thesis.\footnote{The name  {\it h-principle} was chosen a few years later.} Probably, Thurston had noticed a precursor  of this principle in Smale's 1957 announcement which includes the {\it sphere eversion}. Very likely he also read 
 the written version of Thom's lecture at the Bourbaki seminar on this topic (December 1957).
 That report contained the very first figure illustrating immersion theory; this was a {\it corrugation}.
 
 We remember Thurston explaining to us on a napkin in a Parisian bistro how to create an immersed curve in the 
 plane out of a singular plane curve equipped with a non-vanishing vector field: just make corrugations (a kind of waves) along the curve
 in  the direction of the vector field. The beautiful 1992 pamphlet by Silvio Levy, ``Making waves, A guide to the ideas
 behind \emph{Outside In}'', contains an expository paper by Thurston on corrugations with application to the 
 classification of immersed plane curves (the Whitney-Graustein Theorem) and above all, a few steps of 
  a sphere eversion.
  
  Vincent Borrelli (from Claude-Bernard University in Lyon) applied the idea of corrugation in a geometric context
to the problem of finding isometric embeddings in the  $C^1$ category.  After the work of Nash, as generalized by Kuiper, this problem had a 
  theoretical solution: such isometric embeddings exist. However there was no concrete method of constructing them.  
  Borrelli used corrugations together with the so-called {\it convex integration} method of Gromov to find an algorithm consisting of a succession of corrugations and convex integrations for building a $C^1$-embedding of the flat torus into 3-space. 
 Unlike the Nash--Kuiper existence result, Borrelli's algorithm can be implemented on (big!) computers and produces pictures of such a flat torus. This was  done in  collaboration with computer scientists.
   
\begin{figure}[htbp]
\centering
\includegraphics[scale=.12]{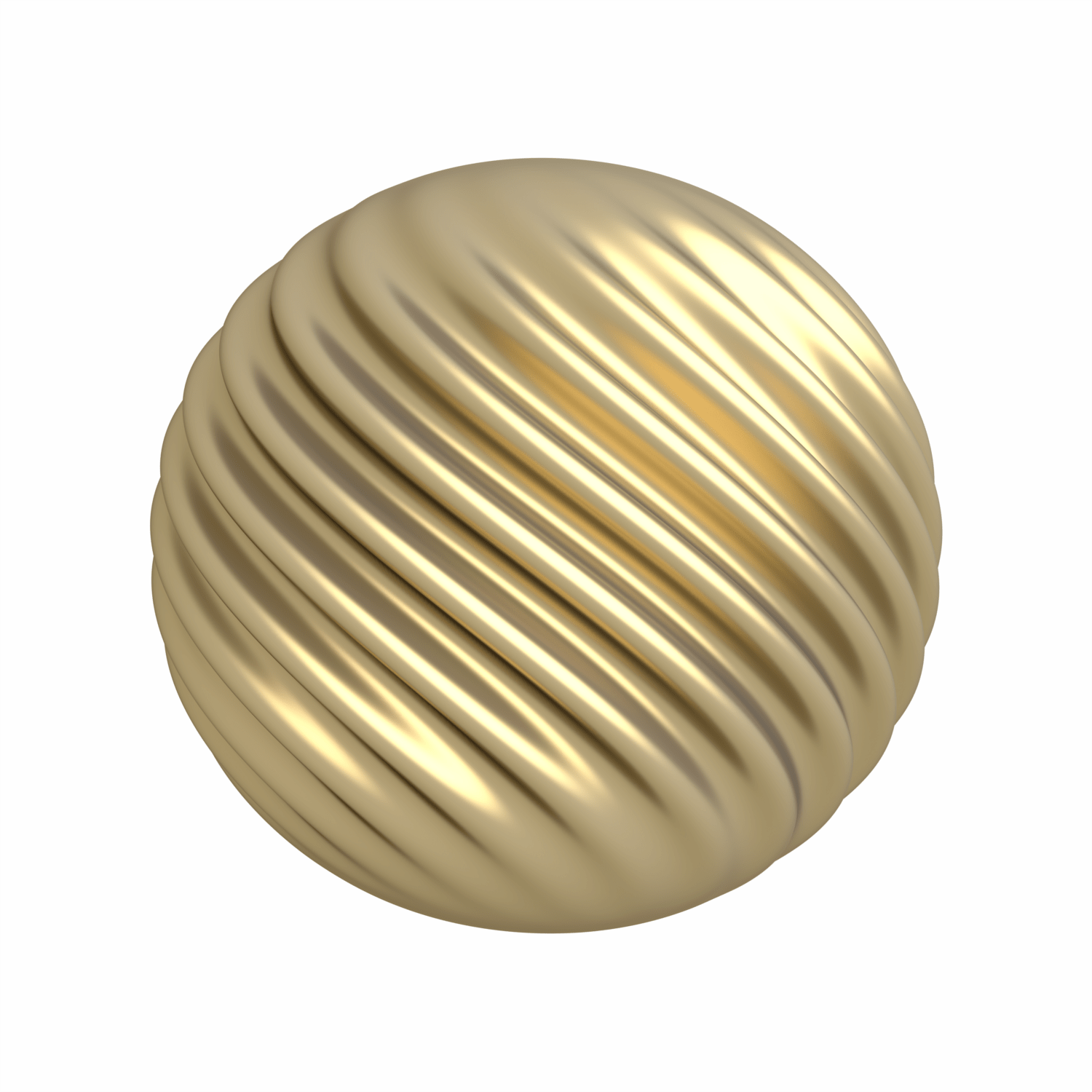}
\caption{First corrugating step for an isometric embedding of the unit sphere into the ball of radius 1/2.
Courtesy of the Hevea Project.}\label{Hevea}
\end{figure}

   After this initial success, Borrelli obtained a $C^1$-embedding of the unit sphere into a ball of radius 1/2 (see Figure \ref{Hevea}).   In a recent paper, his student M\'elanie Theilliere considerably simplified the convex integration method
 to obtain a ``lighter'' algorithm. 
   
              \section{One-dimensional dynamics}

 There is a topic in dynamics which we still have not talked about, that Thurston started working on with Milnor in the year (1972--1973) he spent  at the Institute for Advanced Study, namely, 1-dimensional dynamics. Thurston and Milnor  studied the quadratic family of maps $x\mapsto x^2+c$ and  showed that it is universal in the sense that it  displays the dynamical properties of any unimodal map.\footnote{A paper by Milnor and Thurston on this subject was published in 1988 \cite{MT}. Leo Jonker, in reviewing this paper in Mathscinet writes: ``If there were a prize for the paper most widely circulated and cited before its publication, this would surely be a strong contender. An early handwritten version of parts of it was in the reviewer's possession as long ago as 1977."}  Several important concepts emerged from this work, including the notion of universality, kneading sequence, and kneading determinant.   This work became one of the themes of the seminars that Sullivan conducted at IH\'ES, where he became a permanent member in  1974. 
 
 Sullivan's seminar had a considerable impact on French
  mathematical physicists, and activity on this topic  continues today. In a paper published in 1996, 
 Viviane Baladi (Paris-Sorbonne) and David Ruelle (IH\'ES) revisited the Milnor--Thurston
determinants (as they are called today) in more general 1-dimensional settings  \cite{BR}. These give the entropy of a piecewise monotone map of the interval in terms of the smallest zero of an analytic function, keeping track of the relative positions of the forward orbits of critical points. Baladi worked out a generalization in higher dimensions, which is the subject of her recent book \emph{Dynamical zeta functions and dynamical determinants for hyperbolic maps. A functional approach} \cite{B}. 

Among the other works in France  that further develop Thurston's work on 1-dimensional dynamics, we mention that of Tan Lei and Hans Henrik Rugh (Orsay), \emph{Kneading with weights}, in which  they generalize  Milnor-Thurston's kneading theory to the setting of piecewise continuous and monotone interval maps with weights \cite{RL}.    Rugh's paper \emph{The Milnor-Thurston determinant and the Ruelle transfer operator}  \cite{Rugh} gives a new point of view on the Milnor--Thurston determinant.
  
  By the end of his life,  Thurston returned to the study of the dynamics of maps of the interval. This was the topic of his last paper (published posthumously) \emph{Entropy  in dimension one} \cite{Thurston2014}. The paper uses techniques from a number of other fields on which he worked:   train tracks, zippers, automorphisms of free groups, PL and Lipschitz maps, postcritically finite maps, mapping class groups and a generalization of the notion of pseudo-Anosov mapping, Perron-Frobenius matrices, and Pisot and Salem numbers (two classes of numbers that appear in Thurston's theory of pseudo-Anosov mapping classes of surfaces).

   \section{Geometric structures}

        Foliations and contact structures are closely related to the notion of \emph{locally homogeneous geometric structures} introduced by Ehresmann (in the paper \cite{Goldman}, these structures are called Ehresmann structures). This notion, with its associated developing map and holonomy homomorphism  constitutes one of the key ideas revived by Thurston and that we find throughout his work on low-dimensional manifolds. It originates in the work of Ehresmann from 1935, which is based on earlier contributions by
\'Elie Cartan and Henri Poincaré.  We refer the interested reader to the recent geometrico-historical article by Goldman \cite{Goldman}. Thurston framed the Geometrization Conjecture in the context of locally homogeneous geometric structures, thereby rejuvenating interest in this field of mathematics.
He also developed the theory of geometric structures with orbifold singularities.

    Singular flat structures on surfaces with conical singularities whose angles are rational multiples of right angles provide examples of orbifold geometric structures. Thurston had been developing that theory, including the relation with interval exchange transformations, billiards and Teichm\"uller spaces, since his early work on surfaces, in the mid-seventies, and he motivated the work of Veech and others on the subject. 
    
    In a preprint circulated in 1987 titled \emph{Shapes of polyhedra},\footnote{The paper was published in 1998 under the title \emph{Shapes of polyhedra and triangulations of the sphere} \cite{Thurston-shapes}.} Thurston studied moduli spaces of singular flat structures on the sphere, establishing relations between these spaces and number theory. At  the same time, he introduced the notion of $(X,G)$-cone-structure on a space, extending to the orbifold case the notion of $(X,G)$-structure. In particular he proved that the space of Euclidean cone structures on the 2-sphere with $n$ cone points of fixed cone angles less than $2\pi$ and with area 1 has a natural K\"ahler metric which makes this space  locally isometric to the complex hyperbolic space $\mathbb{C}H^{n-1}$.
     He also proved that the metric completion of that space is itself a hyperbolic cone manifold. This work provides geometric versions of  results of Picard, Deligne and Mostow on discrete subgroups of $\mathrm{PU}(n,1)$, interpreting them in terms of flat structures with conical singularities on the sphere.  Over the years since, the subject of flat structures on surfaces with conical singularities and their moduli  has become an active field of research in France (there are works by Pascal Hubert, Erwan Lanneau, Samuel Leli\`evre, Arthur Avila and many others.)

Influenced by Thurston's ideas, the study of foliations with transverse geometric structures emerged as another research topic among geometers in France. One should add that in the 1970s, and independently of Thurston's work, several PhD dissertations on this topic were defended in Strasbourg, under the guidance of Reeb and Godbillon, including those of Edmond Fedida on Lie foliations (1973), Bobo Seke on transversely affine foliations (1977), and Slaheddine Chihi on transversely homographic foliations (1979).   After Thurston gave several new interesting examples of transverse structures of foliations (the typical one is the class of singular foliations of surfaces equipped with transverse measures)  the theory became much more widely studied, and some geometers started working on foliations equipped with a variety of transverse structures; these included Isabelle Liousse on transversely affine foliations, Ga\"el Meigniez and Thierry Barbot on transversely projective foliations,  \'Etienne Ghys and Aziz El Kacimi-Alaoui on transversely holomorphic foliations and  Yves Carri\`ere on transversely Riemannian and transversely Lie foliations.  Several other mathematicians (such as Abdelghani Zeghib and Cyril Lecuire) began working on laminations in  various settings. The notion of complex surface lamination also emerged from Thurston's ideas and  was   studied  by   Ghys, Bertrand Deroin, Fran\c cois Labourie and others.

 Thurston was also the first to highlight the importance of the representation variety $\mathrm{Hom} (\pi_1(S),G)$, where $S$ is a surface and $G$ a Lie group, in the setting of geometric structures, a point of view which eventually gave rise to the growing activity on higher Teichm\"uller theory. He was the first to realize explicitly that holonomy of geometric structures provides a map
from the deformation space of Ehresmann structures into the representation variety,
which tries to be a local homeomorphism. Although many examples in specific cases of this were known
previously, Thurston realized that this was a very general guiding principle for the classification
of locally homogeneous geometric structures.  We refer the reader to Goldman's article \cite{Goldman} in which he talks about what he calls the \emph{Ehresmann-Weil-Thurston holonomy principle}. Labourie, McShane, Vlamis and Yarmola  in their papers used the expression ``Higher  Teichm\"uller-Thurston theory", and this is likely to become a generally accepted name.

 From a philosophical point of view, Thurston was an intuitionist, a constructive and an experimental mathematician.   He was also among the first to use computers in geometry, in combinatorial group theory and in other topics, and to talk about the rapidity of convergence of geometric construction algorithms.  During a visit to Orsay in November 1987, he gave three talks in which computing played a central role.\footnote{The titles were \emph{Automatic groups with applications to the braid group}, \emph{Conway's tilings and graphs of groups} and \emph{Shapes of polyhedra}.} The book \emph{Word processing in groups} \cite{WPG}, written by Cannon, Epstein, Holt, Levy, Peterson and Thurston, is the result of Thurston's ideas on cellular automata and  automatic groups. These ideas formed the basis of the work of several researchers in France  (Coornaert in Strasbourg, Short and Lustig in Marseille, etc.)
  
  \section{Grothendieck}
 
We cannot speak of Thurston's influence in France without mentioning Alexander Grothendieck, the emblematic figure who worked at IH\'ES for a dozen years and then resigned in 1970 on the pretext that the institute was partially run by military funds. One may note here that Thurston was similarly involved in the US in a campaign against military funding of science. In the 1980s, the \emph{Notices of the AMS} published several letters from him on this matter.   In an attempt to obtain a position at the French CNRS in the years that followed,\footnote{The application was unsuccessful.} Grothendieck wrote his famous research program called \emph{Esquisse d'un programme} \cite{Esquisse} (released in 1984), in which he introduced his theory of \emph{dessins d'enfants} and where he set out the basis for an extensive generalization of Galois theory and for what later became known as Grothendieck--Teichm\"uller theory. At several places of his manuscript Grothendiek expresses  his fascination for Thurston-type geometry, drawing a parallel between his own algebraic constructions in the field $\mathbb{Q}$ of rational numbers and what he calls Thurston's ``hyperbolic geodesic surgery" of a surface by pairs of pants decompositions. He also outlined a principle which today bears the name \emph{Grothendieck reconstruction principle}. This principle
 had already been used (without the name) in the 1980 paper by Hatcher and Thurston \emph{A presentation for the mapping class group of a closed orientable surface} \cite{HT} in the following form:  there is a hierarchical structure on the set of surfaces of negative Euler characteristic ordered by inclusion in which ``generators" are 1-holed tori and 4-holed spheres and ``relators" are 2-holed tori and 5-holed spheres. The analogy between Grothendieck's and Thurston's theories was expanded  in a paper by Feng Luo \emph{Grothendieck's reconstruction principle and 2-dimensional topology and geometry} \cite{Luo}.  Incidentally, the result of Hatcher and Thurston  in the paper  mentioned above is based on Cerf theory.  We mention that Jean Cerf was a professor at Orsay. He was appointed there in the first years of existence of that department, and he created there the topology group (at the request of Jean-Pierre Kahane).
 
 Grothendieck's ideas on the action of the absolute Galois group and on profinite constructions in Teichm\"uller's theory that are based on Thurston-type geometry are also developed in his  \emph{Longue marche \`a travers la th\'eorie de Galois},  a 1600-page manuscript completed in 1981 which is still unpublished \cite{Longue}.  At the university of Montpellier, where he worked for the last 15 years of his career, Grothendieck conducted a seminar on Thurston's theory on surfaces, and directed Yves Ladegaillerie's  PhD thesis on curves on surfaces. 
 
 Grothendieck again mentions Thurston's work on surfaces in his mathematical autobiography, \emph{R\'ecoltes et semailles} \cite[\S 6.1]{RS}.  In that manuscript he singles out twelve themes that dominate his work and which he describes as ``great ideas" (grandes id\'ees).  Among the two themes he considers as being the  most important is what he calls the ``Galois--Teichm\"uller yoga", that is, the topic now called Grothendieck--Teichm\"uller theory \cite[\S 2.8, Note 23]{RS}.
    
      \section{Bourbaki seminars} 
 
Thurston's work was the subject of several reports at the \emph{S\'eminaire Bourbaki}.  This seminar is held three times a year in Paris (over a week-end). It is probably still  the most attended regular mathematical seminar in the world.

In the three academic years 1976/1977, 1978/1979 and 1979/1980, a total of five Bourbaki seminars were dedicated to Thurston's work. In the first one, titled \emph{Construction de feuilletages, d'apr\`es Thurston} \cite{Roussarie}, Roussarie reports on Thurston's result saying that any compact manifold without boundary whose Euler characteristic vanishes admits a $C^\infty$ foliation of codimension one. 
In the second seminar, titled \emph{$B\Gamma $ (d'apr\`es John N. Mather et William Thurston)} \cite{Sergeraert}, Sergeraert reviews one   
of Thurston's deep theorems: 
the homology of the group of $C^r$ diffeomorphisms of $\mathbb{R}^q$ with compact support (as a discrete group)
is closely related to the homology of $\Omega^q(B\Gamma^q)$, the $q$-th loop space of the classifying space of
codimension-$q$ Haefliger structures of class $C^r$.\footnote{The notion of Haefliger structures translates an idea of 
singular foliation equipped with a desingularization. This leads to a homotopy functor which has a classifying space, 
analogous to $BO(n)$ for vector bundles of rank $n$. The case $q=1$ was already known to J. Mather in 1970. For this reason, 
one speaks today of the Mather-Thurston homology equivalence theorem.
 Let us mention that Takashi Tsuboi circulated a pamphlet with pictures
explaining a map that induces this isomorphism. It is available on Tsuboi's homepage.} In the third seminar, titled \emph{Travaux de Thurston sur les diff\'eomorphismes des surfaces et l'espace de Teichm\"uller} \cite{Po}, Po\'enaru gives an outline of Thurston's theory of surfaces, which appeared later in \cite{FLP}. In the fourth seminar, titled \emph{Hyperbolic manifolds (according to Thurston and J\o rgensen)} \cite{Gromov},  Gromov reports on some of the powerful techniques contained in Thurston's  1997/98 Princeton notes, including his work on limits of hyperbolic 3-manifolds, his rigidity theorems, and the result stating that the set of values of volumes of hyperbolic 3-manifolds of finite volume is a closed non-discrete subset of the real line. As a matter of fact, Gromov arrived to France and lectured at Orsay at the end of the 1970s. His notion of simplicial volume played, via the techniques of smearing out and straightening, a key role in the (so-called Gromov--Thurston) version of Mostow rigidity theorem for 3-dimensional hyperbolic manifolds contained in Chapter 6 of Thurston's Princeton notes \cite{Thurston1}.
In the fifth Bourbaki seminar, titled \emph{Travaux de Thurston sur les groupes quasi-fuchsiens et les vari\'et\'es hyperboliques de dimension 3 fibr\'ees sur $S^1$} \cite{Sullivan}, Sullivan gives an outline of Thurston's results on hyperbolic structures on irreducible 3-manifolds which fiber over the circle and which contain no essential tori.  At the same time, using a limiting procedure in the space of quasi-Fuchsian groups, he gives a new proof of Thurston's result saying that the mapping torus of a homeomorphism of a closed surface of genus $\geq 2$ with pseudo-Anosov monodromy carries a hyperbolic structure.

Thurston's work has been the subject of several other Bourbaki seminars over the years;   we mention in particular seminars by Morgan on finite group actions on the sphere \cite{Morgan}, by Ghys on the Godbillon-Vey invariant \cite{Ghys}, by Boileau \cite{Boileau} on uniformization in dimension three, by Lecuire on ending laminations \cite{Lecuire}, and, finally, by Besson on the proof of the geometrization  of 3-manifolds  and the Poincar\'e conjecture \cite{Besson}. 
 
\section{The last visits to Paris}

Speaking of the proof of the Poincar\'e conjecture---another problem rooted in French mathematics, a problem that haunted Thurston during all his mathematical life--- we are led to the last time we saw Thurston 
 in Paris. This was in June 2010 at the Clay research conference, where he gave two beautiful talks at the magnificent lecture hall of the  Oceanographic Institute.  The first talk was titled ``The mystery of three-manifolds." The second one, shorter, was a \emph{Laudatio} on Grigory Perelman. Thurston recounted his personal experience with the Poincar\'e conjecture. In a few minutes, he expressed his deep admiration and appreciation for Perelman and he said in a few moving words how much he was gratified to see that the geometrization conjecture became a reality during his lifetime. With an amazing humbleness, he declared that when he read the proof he realized that it is a proof that he could not have done (``some of Perelman's strengths are my weaknesses"). He concluded with these words: 
\begin{quote}\small
Perelman's aversion to public spectacle and to riches is mystifying to many. I have not talked to him about it and I can certainly not speak for him, but I want to say I have complete empathy and admiration for his inner strength and clarity, to be able to know and hold true to himself. Our true needs are deeper---yet in our modern society most of us reflexively and relentlessly pursue wealth, consumer goods and admiration. We have learned from Perelman's mathematics. Perhaps we should also pause to reflect on ourselves and learn from Perelman's attitude towards life.
\end{quote} 
\begin{figure}[htbp]
\centering
\includegraphics[scale=.3900]{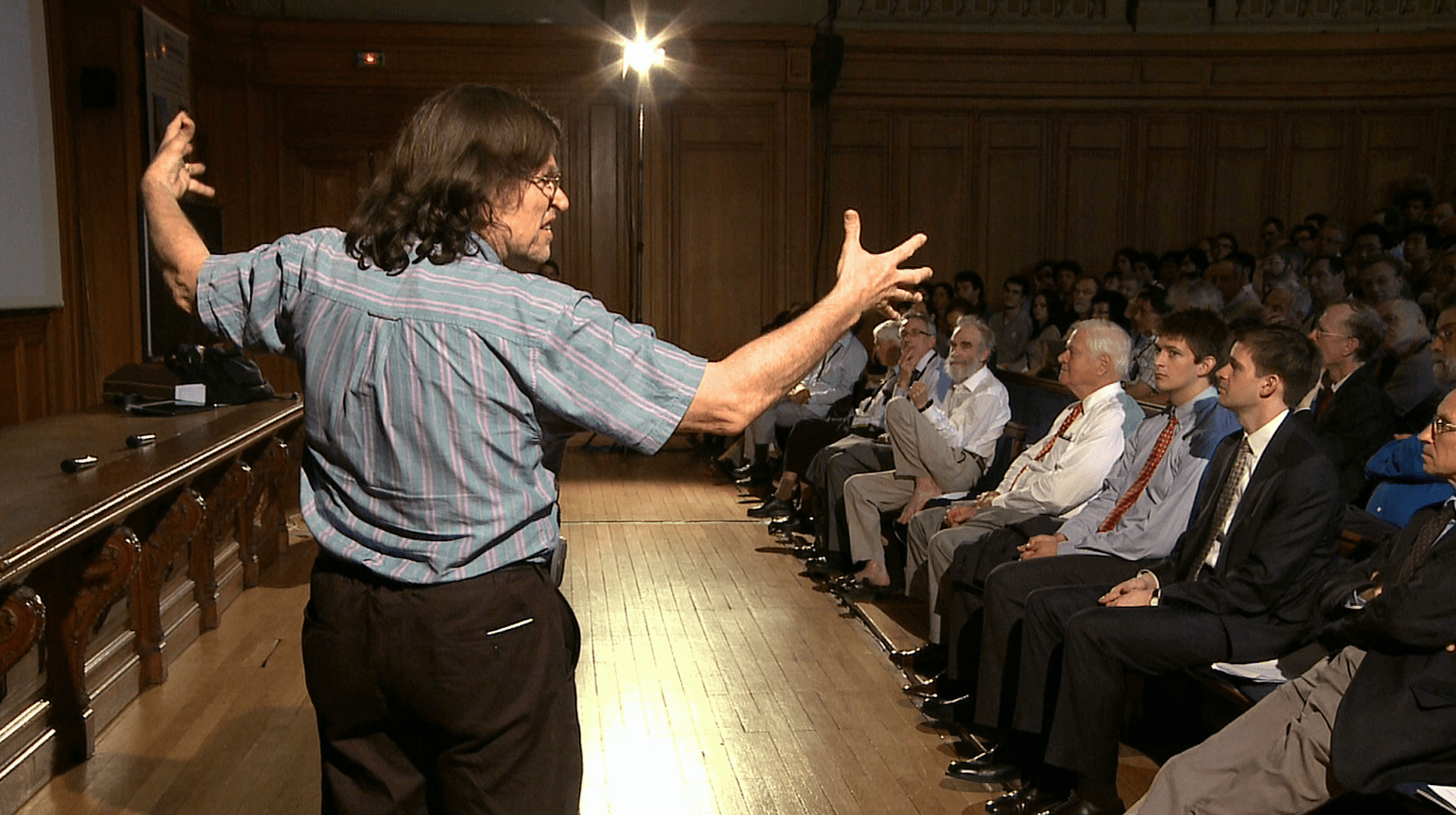}
\includegraphics[scale=.4650]{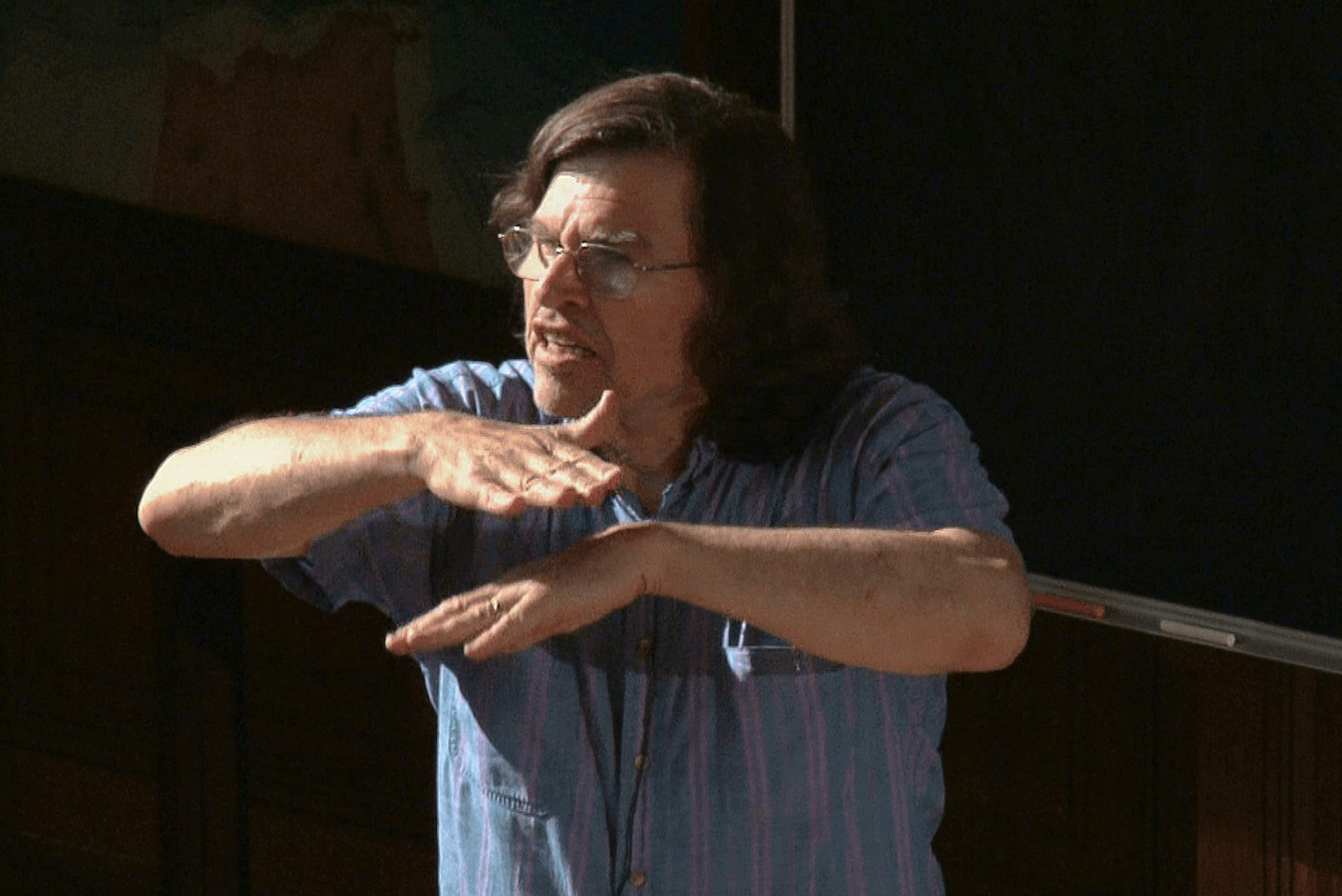}
\caption{Bill Thurston at the Clay conference in Paris, Oceanographic Institute, June 2010. Photos $\copyright$~Atelier EcoutezVoir.}\label{Thurston-Mystery}
\end{figure}

   Paris is also the capital of fashion design.  A few weeks before the Clay conference Thurston was there for  the
   fashion week, which takes place every year in March at the Carrousel du Louvre, which sits between the Jardin des Tuileries and the Louvre museum. 
  The fashion designer Dai Fujiwara presented  a beautiful collection of pieces made for the Issey Miyake brand, inspired by Thurston's eight geometries. 
 A journalist covering the event wrote: ``Two decades ago, in the same venue, Romeo Gigli transfixed Paris with a show so rich and romantic that it moved its audience to tears. Maybe that didn't happen today, but, at the very least, Fujiwara used his inspiration to blend art and science in a manner so rich and romantic, it stirred the emotions in a way that reminded us of Gigli."\footnote{Tim Blanks, \emph{Vogue}, March 5, 2010.} Among the other comments on this event, we note:  ``Fashion scaled the ivory tower at Miyake, where complicated mathematical theorems found expression in fabric" (Associated Press); ``Mathematics and fashion would seem to be worlds apart, but not so, says Dai Fujiwara" (The Independent); 
``you did not need a top grade in maths to understand the fundamentals of this thought-provoking Issey Miyake show: clean geometric lines with imaginative embellishment" (International Herald Tribune); ``Fujiwara used his inspiration to blend art and science in a manner so rich and romantic, it stirred the emotions" (Style).

  Thurston wrote a brief essay, distributed during the Miyake fashion show, on beauty, mathematics and creativity.  
  Here is an excerpt:
 \begin{quote}\small
Many people think of mathematics as austere and self-contained. To the contrary, mathematics is a very rich and very human subject, an art that enables us to see and understand deep interconnections in the world. The best mathematics uses the whole mind, embraces human sensibility, and is not at all limited to the small portion of our brains that calculates and manipulates with symbols. Through pursuing beauty we find truth, and where we find truth, we discover incredible beauty.
  \end{quote}

 In another article written on that occasion  for the fashion magazine \emph{Idom\'en\'ee}, Thurston made the following comment about the collection: 
\begin{quote}\small
The design team took these drawings as their starting theme and developed from there with their own vision and imagination. Of course it would have been foolish to attempt to literally illustrate the mathematical theory--- in this setting, it's neither possible nor desirable. What they attempted was to capture the underlying spirit and beauty. All I can say is that it resonated with me.
  \end{quote}
  In an interview released on that occasion Thurston recounted how he came to contribute to the collection, and he declares there: ``Mathematics and design are both expressions of human creative spirit.''
  One of the comments on the video posted after this interview says: ``I can't believe this mathematics guy. He's so ... not like what I expected.''

 \part{}
 
   \section{Valentin Po\'enaru}

   Bill Thurston went through the mathematical sky like an immensely bright, shiny, meteor. My short contribution certainly does not have the ambition of trying to describe his fantastic trajectory. Rather, more modestly, we want to tell
something about Bill's impact on the mathematical life of the Orsay Department of Mathematics (Universit\'e de Paris-Sud) in the seventies.
   
We heard for the first time about Thurston's activity from our friend
and then colleague Harold Rosenberg, who invited Bill to talk about his work
on foliations. The first contact already triggered the mathematical career of 
Michel Herman.
   
What we had heard so far was Bill's work on foliations of codimension higher than one. Next, we still vividly remember the lectures which Bob Edwards gave us on Thurston's codimension-one theorem. Several other mathematicians had
already tried to explain this to us, but they always got bogged down before managing
to get to the main point. Bob's lectures were not only very illuminating,
they also had the aesthetic quality of a drama stage, with suspense and ``coup de th\'e\^atre."
   
Then, via Dennis Sullivan who was visiting Orsay for a year, before settling
for the next twenty years or so at the neighboring IH\'ES, we were introduced
to the hyperbolic world, where Bill was now one of the brightest stars. We 
vividly remember how this field where Bill, Misha Gromov, Dennis and others did such big things, was then looked down upon by many. A distinguished colleague told us that hyperbolic geometry was just ``a gadget." In his
panoramic books on the mathematics of our time, Dieudonn\'e does not even
mention the topic, obviously thinking that it was too marginal and parochial.
   
It is also Dennis who, via some very convincing drawings, conveyed to us
Thurston's discovery that, by infinite iteration via a pseudo-Anosov diffeomorphism,
curves can turn into measured foliations. And when we tried to make
that quantitatively precise, a very nice mixing property popped up; ergodic
theory was there, big.
   
Somebody brought us some notes written after some lectures of Thurston on his theory of surfaces. This triggered the
Orsay seminar on this topic, organized by Fathi, Laudenbach and myself, where we tried to provide full
details for the theory. A lot of distinguished visitors joined us. But then,
there was a big problem where we got bogged down. We did not manage to
glue the space of measured foliations to the infinity of the Teichm\"uller space.
The Teichm\"uller specialists offered us several suggestions, but none of them was
good, since we needed a \emph{natural} compactification of the Teichm\"uller space where the automorphisms of the surface should extend continuously.

   Then, during a very hectic week at Plans-sur-Bex, in the Swiss Alps, we were together with Bill who during a memorable and intense half-hour gave
us the correct hint on how to proceed. Thus, we could both finish our seminar
and write the corresponding book. The seminar in question, and in
particular the influence of Albert Fathi, triggered the mathematical activity of
Jean-Christophe Yoccoz.
   
Next, we were introduced by Sullivan, and others, to Thurston's program
of introducing hyperbolic geometry into the field of 3-manifolds. With
Larry Siebenmann as the main organizer, this created a lot of activity. This is
how Francis Bonahon and Jean-Pierre Otal started their mathematical trajectories. Our colleagues and friends Fran\c cois Labourie and Pierre Pansu were also
quite influenced by all these activities.
   
Very naturally, our Department of Mathematics proposed to our University that
Bill be awarded the title of Doctor Honoris Causa. The ceremony was held in the
Fall of 1986. It so happened that at about the same time, the Computer Science
Department did the same for Donald Knuth. Thus, in the largest lecture room
on our campus there were two big public lectures by the two laureates. Bill gave
us a brilliant and amazing lecture on how geometry (hyperbolic or SOL) can
optimize computer construction, providing maximum connectivity with a minimum of
parasitical interferences, in a given volume.
   
In the beginning of the Summer of 2010 Bill came for a last time to Paris
and lectured on some future plans of his at the Institut Henri Poincar\'e. But
those were not to be since Bill passed away soon afterwards.
   
Bill Thurston's impact on the mathematical life of our department was immense
and quite a number of people here owe a lot to him. We will never forget
him.

 \section{Harold Rosenberg}
 
 (Excerpt from an email dated March 22, 2016, addressed to Fran\c cois Laudenbach)

  In 1971, I visited Berkeley for 6 months.  I taught a class on
  foliations, and Bill Thurston was a student in the class.
  By the end of the term, Bill and I were working together, and we wrote a
  paper which was  published in the proceedings of a meeting on dynamical systems  which took place in the city of Salvador, on July 26-August 14, 1971. (This
 was my first visit to Brazil.) The book was edited by M.
  Peixoto.
  Upon returning to Paris, I invited Bill to visit.  He came to Orsay.
  Haefliger organized a meeting (I believe it was in Plans-sur-Bex), and
  Bill and I went and Haefliger met Bill.  Also, Bill met Dennis in Paris at that time.
  Bill was working on the geometry of foliations in 3-manifolds at the
  time.  Shortly thereafter, he managed to construct a family of
  foliations whose Godbillon-Vey invariants varied continously
  onto an interval.  Based on this, Milnor offered Bill a visiting position
  at Princeton.  The following year (Bill was in Princeton), he started
  proving his important integrability
  theorems in higher dimensions, and thinking about diffeomorphisms of
  surfaces. 
  
    \section{Francis Sergeraert} 
  (Excerpt from an email dated April 4, 2016, addressed to Fran\c cois Laudenbach, translated from the French)

I remember very well Thurston's visit to Orsay in 1971. It was the time of his proof of $\mathrm{Im(GV)} = \mathbb{R}$. Someone (I don't remember who, maybe Rosenberg) was trying to find a mistake in that result.  

At that time, Thurston already knew how to connect the homological aspect of diffeomorphism groups and that of 
  $B\Gamma$. Michel Herman's note on the simplicity of  $\mathrm{Diff}_+(T^n)$ dates from the same year. I don't know what was the real influence of Thurston on that work, but it was non-negligible. The spaces  $B\Gamma$, since their discovery by Haefliger, were still very mysterious. Mather had started to unblock the problem by first establishing the connection with the homology of $\mathrm{Diff}_c(\mathbb{R})$. He had difficulties publishing his result, which eventually appeared in 1975, but 
  starting in 1973, Thurston, in Princeton, was already explaining its generalization in any dimension.
   
   Regarding the Poincar\'e conjecture, when Thurston started working on it, everybody considered it as a problem in algebraic topology. Thurston's work on foliations was probably the first indication that there were connections with differential analysis. It is tempting to see there the germ of the idea of geometrization.

     \section{Norbert A'Campo} 
 
 (Excerpt from an email dated February 1, 2019, addressed to Athanase Papadopoulos, translated from the French)
 
 In the beginning of the seventies, an offer was made to Morris Hirsch for a position at Orsay. He accepted, paperwork was done, but Hirsch eventually resigned, the reason being that he had a very good student which he did not want to leave. After a moment of consternation, Rosenberg and Siebenmann decided to invite the student. Thurston came to Orsay.  He also visited Dijon and Plans-sur-Bex. We did a large portion of these trips was together, by car.

The conference at Plans-sur-Bex was organized by Haefliger, one week in March.
The wonderful inspiring  wine was provided by Kervaire. Local organisation with excellent cooking was done by the family Amiguet from Geneva.

During the travel, and at the conference, Thurston explained during his talks and in informal discussions some vivid and beautiful new mathematics. In particular he proved that the Godbillon-Vey map
$$GV: \{\textit{codimension 1 foliations on}\, S^3\}
 \to \mathbb{R}$$ is surjective. He used the theorem stating that planar hyperbolic polygons of equal area are scissor equivalent. For me, this was the first time I saw hyperbolic geometry at work.

A few years later, after many other celebrated uses of and contributions to hyperbolic geometry by Thurson, Larry Siebenmann asked me to give a graduate course on
hyperbolic geometry at Orsay. I knew nothing about that subject. Fortunately, I
was planning a visit the Mittag-Leffler Institute, and in one of the attics there I
found several old documents on this topic. I took some notes and I came back to
Orsay ready to give my course and (at the same time) learn hyperbolic geometry.

  \section{Gilbert Levitt}
 
I first heard the name Thurston in 1976 (or maybe 1975). David Epstein gave a graduate course on foliations in Orsay. He explained Thurston's stability theorem for foliations of codimension 1, and a large part of the course was devoted to trying to understand the proof of Thurston's theorem on foliations of codimension greater than 1.

The following year I started working towards a PhD under the guidance of Harold Rosenberg. He made me study Thurston's thesis, about foliations of 3-manifolds which are circle bundles, and my first published paper may be viewed as a write-up of that thesis. 

Following the advice of Harold,    I then turned to (singular) foliations on surfaces. Almost all  my papers on that subject mention Thurston in the bibliography; when compactifying Teichm\"uller space and classifying homeomorphisms of surfaces, he introduced measured foliations (and laminations), which may be viewed as a building block for constructing general foliations on surfaces. 

He also defined train tracks, which I encountered later while working on automorphisms of free groups. Following the seminal paper by Bestvina-Handel, train tracks were carried over from surfaces to free groups and became an extremely important tool in geometric group theory.

One of Thurston's last contributions is a paper posted on the arxiv in 2014, where (among other things) he completely characterizes which   numbers may appear as growth rates of automorphisms of free groups.

My first meeting with Thurston was in the fall of 1978. At that time Harold Rosenberg visited Santa Cruz for several months and he arranged for R\'emi Langevin and me to spend some time in Berkeley. On the way home we stopped in New York and he secured an appointment with Thurston for me. 

R\'emi and I drove out to Princeton and I spent quite some time with Thurston discussing foliations on surfaces and related topics. For some reason we didn't use a blackboard, and I still have a notepad covered with his drawings. I was 23 at the time, I hadn't yet proved a real theorem, and thinking back I am really grateful and honored that he devoted so much of his time to me.

  \section{Vlad Sergiescu}

I would like to say a few words on Thurston's influence on French foliation theorists in (and between) Lille and Orsay, in a specific situation: the geometry of the Godbillon-Vey class of a codimension-one foliation.

Around 1970 a big excitement arose in foliation theory, due to the discoveries of the Bott vanishing theorem, Haefliger's classifying space, Fuchs-Gelfand cohomology and the Godbillon-Vey class.

For a foliation given by a 1-form $\omega$, let $\eta$ be a 1-form such that $d\omega=\omega\wedge \eta$. The 3-cohomology class of the closed form $\eta\wedge d\eta$ is the Godbillon-Vey class $GV$.

A short time after its discovery, Harold Rosenberg, who introduced foliations at Orsay and was an outstanding advisor there, wrote a paper with Thurston, published in the proceedings of a conference in Bahia, in which they asked whether the Godbillon-Vey class vanishes for a foliation without holonomy.

In 1973, Sullivan began a course at Orsay on his new rational homotopy of differential forms, with (among others) the following question: What does the Godbillon-Vey class measure?

Many people around Orsay contributed in some way or another to make advances on this question. 
 Let me mention Roussarie, Herman, Moussu, Sergeraert, and Roger as well as Haefliger, Epstein and Sullivan as regular visitors.  Roussarie proved that the GV class is non-zero on the unit tangent bundle of a hyperbolic surface. Thurston showed this as well (as a graduate student), and later, he proved that it varies continuously! In 
 an influential
paper, Herman proved (and Guy Wallet as well) a vanishing theorem on the torus $T^3$.

 During those years, Thurston was involved in foliations, until 1976--77, when he switched to 3-manifolds where  
 he made a strong use of measured foliations for diffeomorphisms of surfaces. He revolutionized both subjects. Some of his landmark theorems, besides the continuous variation, are the existence of a codimension-one foliation on any manifold with vanishing Euler characteristic and the $(q+1)$-connectedness of the classifying space $B\Gamma_{q}$. Let me add his inspiring picture of the helical wobble.

 In 1976, I joined a seminar in Lille organized by Gilbert Hector. Among the
participants  were Duminy, Ghys, El-Kacimi Alaoui, Lehmann, and soon, first and
foremost, Alberto Verjovski. All of us were sensible to the GV world. There, we learned that  the absence of Lamoureux' resilient leaves (which are self-spiraling) is a good context to attack the problem of the vanishing of GV. It contained the non-exponential setting that Moussu and Sullivan had already suggested (thus also the vanishing of the holonomy). 
  
Shortly after a joint note for the $S^1$-foliated case, G\'erard Duminy found a brilliant proof of the general vanishing theorem. One major innovation was the exploration of a decomposition of GV as the product of a ``Godbillon" measure and a ``Vey" class. Several other new ideas paved the way for the introduction of techniques  of ergodic theory and dynamical systems in foliation theory. An example is the connection between the entropy defined by Ghys--Langevin--Walczak and the GV class.
Hurder and Katok proved that GV is invariant under absolutely continuous homeomorphisms, while the
topological invariance is still open---Ghys obtained several results in this direction. Hurder and Langevin present a
modern view of the above topics in a recent article on GV and $C^1$
dynamics.

  Thurston introduced a 2-cocycle $t_{gv}$ on $\mathrm{Diff}(S^1)$ (the names of Bott or Virasoro are sometimes linked to this) closely related to GV:
  
   \[(f,g)\to \int \log g'(\log(f\circ g)')'.  \]

A useful observation made later was that $t_{gv}$ can be extended below the $C^2$ class, to the Denjoy $P$ class (maps with bounded
log-derivative variation)\footnote{This terminology was used by M. Herman.} and to the class $C^{1+\alpha}$, $\alpha>\frac{1}{2}$. 
Takashi Tsuboi made a thorough study of such extensions leading to his beautiful GV-cobordism characterization.

 In his Berkeley years, Thurston met a fellow student, Richard Thompson who was working in algebraic logic and discovered in this context three groups with wonderful properties. One of them turned out to be isomorphic to the group T of PL dyadic homeomorphisms of the circle.
 
 Thurston talked about this to several people around him.  Ghys and myself learned about $T$ from Epstein and Sullivan. We then found a PL version of Thurston's cocycle called $\overline{t_{gv}}$. It is intriguing and remarkable that together with the Euler class, it generates the cohomology of $T$, and this turned out to be exactly the Fuchs-Gelfand cohomology of $\mathrm{Diff}(S^1)$.
 
 To conclude, let me point out that the classes $GV$, $t_{gv}$ and $\overline{t_{gv}}$ appear to be ubiquitous, with connections to braids, mapping class groups and Teichm\"uller spaces, loop groups, Chern-Simons invariants, Virasoro algebra and groups, $\mathbb{C}^*$-algebras, index theory, strings, solitons and hydrodynamics.
 
 I never had the privilege to meet Thurston, but my friend and collaborator Peter Greenberg did. He told me once that an important thing he learned from him was how to play with mathematics. I vividly remember a day of 1987 in Mexico when both of us were
playing with David Epstein to connect Thompson groups with braid groups. At that time, it was not a success. Suddenly David exclaimed: I am sure that Thurston would find it on the spot!

\section{Michel Boileau}

 (Excerpt from an email dated February 20, 2017, addressed to Fran\c cois Laudenbach)

The first time I had the occasion to listen to a lecture of Thurston was at the conference held in Bangor (G. B.) 
in July 1979. He gave four lectures. The first three were about the geometrization conjecture of 3-manifolds and 
its proof for Haken manifolds. The last lecture was about the Smith conjecture whose recent 
proof relied upon the geometrization of Haken manifolds. Thurston motivated his conjecture 
by the fact that it dealt with all 3-manifolds. One could hope that it would be solved within the next thirty years; 
history showed that this was right.

 The second occasion on which I followed lectures by Thurston was at the 1984 Durham conference. Thurston gave a 
 series of lectures on the geometrization of orbifolds, in particular on hyperbolic conical structures and their
 geometric limits. 
 Again, objects and methods presented  in dimension three were completely new.  It was only thirteen years later, when I started with Joan Porti and Bernhard Leeb to write a complete
  proof of this theorem, that I understood
 the ideas that Thurston had tried to communicate in these lectures.

 In the Fall of 1986, William Thurston was awarded the degree of Doctor Honoris Causa from the University of Paris-Sud
 (Orsay). On that occasion, he gave two lectures. In the first one, for a large audience, he explained how to apply
 methods from hyperbolic geometry to computer science. The second one, more specialized, was 
 on the deformation space
 of polygons in the plane. On that occasion, Gromov, from the audience, challenged him with an objection. Thurston answered with 
 the smile he has always had when he tried to communicate his extraordinary vision of geometry.

 His style and his manner of communicating mathematical ideas, though very generous, have frequently
 raised criticism. In my opinion, they rather reveal his highly demanding commitment to the quality of writing. In a talk he gave at a conference in Tokyo in 1998, concerning the proof of his orbifold theorem, Thurston declared: ``I am reproached for not writing enough  but what I have in mind
 is much more beautiful than what I am able to put on paper.''
 
At the same conference, I had the opportunity to ask him whether every 3-manifold could have a conical 
 hyperbolic structure with angles arbitrarily close to $2\pi$. This seemed {\it a priori} impossible for the 3-sphere
 and many colleagues to whom I had asked the question were of the same opinion. Thurston immediately answered: ``Yes'', and he showed me
 a simulation on his computer, precisely for the 3-sphere. This result was proved later by Juan Souto. 
 
\section{Pierre Arnoux}

I started as a mathematics graduate student at Orsay at 1979, by following the lectures of Michael Shub on dynamical systems.\footnote{This was called a DEA (Dipl\^ome d'\'Etudes Approfondies) course, usually attended by graduate students, the year before they choose a subject for their PhD dissertation. But the courses were also sometimes followed by confirmed researchers.}, This was two years after the seminar on
 {\it Thurston's work on surfaces}. At that time at Orsay, one was immersed, without even realising it, in a particular mathematical culture.  
 I  became quickly aware of the classification of surface automorphisms, 
 foliations etc., even if I was far from understanding the proofs.
 
 Michel Herman suggested me to work on interval exchange maps for my PhD. thesis. I always thought of them geometrically,
 associated with surface foliations. At that time, very few explicit examples of surface automorphisms
 were known; most of them were related to coverings of automorphisms of the torus.
 We came across a paper by William Veech in which he started the parametrization of holomorphic forms\footnote{The history is quite complicated: Veech was interested in results of ergodic theory; he had already worked on particular types of interval exchange maps given by skew-product of rotations, and found the general idea of induction in a paper of Rauzy, as generalized continued fractions; he then learnt, apparently from Thurston, of the link between interval exchange maps and measured foliations, and this led him to the proof of almost everywhere unique ergodicity. I think his papers have not been fully read, and they still contain a number of unnoticed results.} (hence,  particular strata 
 of the cotangent bundle of Teichm\"uller space) using  interval exchange maps; that gave a way 
 of building {\it self-similar foliations}. I remember a night trip back from England (possibly Durham)
with Jean-Christophe Yoccoz and Albert Fathi during which they built such an example
with a cubic coefficient (I was rather a spectator than anything else). This was the first example of a
 pseudo-Anosov diffeomorphism which does not arise from a torus automorphism. This was also the kind of things
 that I enjoyed: to build explicit and concrete, yet slightly strange objects. Twenty years later, Maki 
 Furukado, a Japanese mathematician, gave me a model of this foliated surface constructed by sewing
rectangles of striped material; this is not a trivial thing to do, because the singularities of the foliation impose heavy constraints. But you can easily see why the suspension has to be on a surface of genus 2. 

A few years later, around 1984, I came across a paper by  G\'erard Rauzy in which he had
  built a self-similar fractal set associated with a substitution whose similarity coefficient is the same as the one of the 
 pseudo-Anosov example. I thought  that this was more than a coincidence.  
 By using the symbolic models associated with the two systems, 
 I was able to show that the interval exchange which had made possible this example was measurably conjugate  to a rotation of the 2-torus  by a continuous map 
 $\mathbb{T}^1\to \mathbb{T}^2$  
  whose image is a Peano curve filling the 2-torus. It followed easily that, by taking suspensions, 
 the pseudo-Anosov at hand was measurably conjugate by a continuous map to a hyperbolic automorphism of the 3-dimensional  torus $\mathbb{T}^3$.
 
In Orsay, people were also familiar with the work of Adler and Weiss  yielding an explicit Markov partition of the hyperbolic automorphisms of $\mathbb{T}^2$ and showing that these automorphisms are classified up to measurable conjugacy by their  largest eigenvalue (their logarithm is the topological entropy, as well as the measure theoretic entropy for the Lebesgue measure). We also knew that any hyperbolic automorphisms of $\mathbb{T}^n $ has a Markov partition. But, these partitions are not at all trivial; for instance, it follows from a very short paper by Rufus
 Bowen that for $n>2$ their boundary is a fractal set. It is indeed difficult to exhibit explicit examples of Markov partitions, except in the case of surface automorphisms, where they are made of explicit rectangles.
The particular example we had found, by giving a measurable conjugacy between an automorphism of a surface with an explicit Markov partition and  a toric automorphism, showed an explicit partition, with fractal boundary of known dimension, for the toral automorphism.
 
 In this theory, there is a basic example easy to understand, that is, $\mathbb{T}^2$ with its Teichm\"uller
 space which is the hyperbolic plane, and its moduli space, which is the classical modular 
 surface equipped with its geodesic flow. There are two simple ways for generalizing that case: taking higher 
 dimensional tori, or taking surfaces of higher genus. Sometimes, one has the feeling that the two ways 
 are but the same: a hyperbolic automorphism of the $n$-dimensional torus $\mathbb{T}^n$ can be often unfolded to an automorphism of a surface  of the corresponding genus, a little like these {\em kirigami}, Japanese paper flowers which unfold when you put them in water.
 
Now I would like to talk about the period 1992--2019 spent in Luminy-Marseille.
 
 I was more and more attracted by  the arithmetic and combinatoric constructions of Rauzy  and I 
 went to work in the laboratory he had founded in Luminy. But I was still interested in the 
 geometric side of these constructions. Thurston had also written a paper, not published but available
 as a xeroxed preprint, about tiles associated with algebraic numbers; this was parallel  to what we were
 doing with substitutions. 
 
 In Luminy, we worked on symbolic sequences with low  (sublinear) complexity, in particular the Sturmian
 sequences which appear in many different settings: dynamics of rotations on the circle, Farey sequences, 
 dynamics of continuous fractions, and more curiously, dynamics of the Mandelbrot set.\footnote{S. Bullett and P. Sentenac wrote a beautiful paper on this subject \cite{Bullett}.} Here, the basic lemma, attributed to Thurston, states the following: 
 an orbit of the map $x\mapsto 2x\  {\rm mod}\ 1$ is cyclically ordered if and only if its binary encoding is 
 Sturmian. In all the papers I have worked on since, there is the influence of the geometry I learnt at Orsay at that time, mixed with the discrete mathematics and the number theory which was the mark of Rauzy. 
 
 In Marseille, there was another group of mathematicians who were studying outer automorphisms of free groups. They included Arnaud Hilion and Martin Lustig, collaborating with Gilbert Levitt who was in Caen. Outer automorphisms of free groups have a lot of analogies with mapping class groups of surfaces. The substitutions that we studied
 in Luminy were simpler cases of these automorphisms, in the same manner as matrices with
 positive coefficients are simpler than general matrices (Perron-Frobenius). The two groups started to collaborate.
 
 In the articles that I write or read today, I continue to feel what happened at Orsay around 1976, 
 with Thurston, Douady, and Hubbard (who is now regularly in Marseille): tiles associated with automorphisms
  of free groups, generalized Teichm\"uller spaces,  explicit conjugations between  bifurcations for families of 
  continuous fractions or for families of quadratic polynomials, etc.

\section{Albert Fathi}

Bill Thurston's impact on French topologists is certainly one of the best influences on the group.
 
I first heard of Bill when I started my Graduate Studies in 1971. At that time he was already a legend for his work on foliations.

I think I first met Bill at the CIME school on Differential Topology in 1976 at Varenna.
Thurston was, along with Andr\'e Haefliger and John Mather, one of the three people delivering the courses.
It is unfortunate that he never delivered the manuscript of his lectures for publication. My most vivid impression
of this meeting was the private explanation by Andr\'e Haefliger on Thurston's beautiful
geometric argument on how to 
obtain that the (connected component of the) diffeomorphism group of a compact manifold
is perfect from the case of the n-dimensional torus that was previously done by Michel Herman.
This was an Aha! moment: how a deep insight in geometry can circumvent the impossible
adaptation to general manifolds of Herman's work on the torus. It used KAM methods and hard implicit 
function theorems in neighborhoods of irrational translations on the torus.

The work of Bill on diffeomorphisms of surfaces led to the monograph that we edited with
Fran\c cois Laudenbach and Valentin Po\'enaru.

It was Valentin Po\'enaru who drew us to this subject. He came one day from IH\'ES with a set 
of hand-written notes that Mike Handel produced while listening to Bill's course in Princeton.
He convinced Fran\c cois and myself to run a seminar on the subject. This seminar took place
in 1976-77 in Universit\'e Paris-Sud (Orsay).

The group of diffeomorphisms of a surface up to isotopy is called the mapping class group
(of the surface). Bill's work essentially produced a ``best'' representative in each element of the
mapping class group.

Valentin Po\'enaru gave us the Grand Tour on the subject in the first lectures
of the seminar. I started to work almost everyday with Fran\c cois
to be able to understand the details. It took us a couple of years to produce a usable manuscript.
We benefited from advice of Francis Sergeraert who served as a referee.
At this time, Bill used measured foliations rather than measured laminations which appeared
after most of our manuscript was finished. This is why measured laminations are not in the monograph. Anyway, I find it very rewarding that 30 years later, it was still found useful
to have an English version of our monograph. The mathematical world is smaller than we think it is: 
one of the two editors of the English version is Dan Margalit who now is my colleague at Georgia Tech.

Bill's work was a  revolution in the old subject  of classification up to isotopy of surface diffeomorphisms. Before him, there was a remarkable work of Nielsen in the 1930s 
which pointed out the elements of finite order of the mapping class group. However,
nobody really realized the existence and irreducibility of what Bill called
pseudo-Anosov diffeomorphisms. Of course, the fact that Anosov diffeomorphisms were, by that time, extensively studied,  in particular, through the properties of the stable and unstable foliations,
is certainly what motivated Bill to introduce these pseudo-Anosov diffeomorphisms.
Obviously, Nielsen could not have benefited from such a knowledge. 
What was also remarkable in Bill's approach was that he also made strong connections with objects, besides pseudo-Anosov maps, that were subjects of intensive research in Dynamical Systems like interval exchange and entropy. For me, who was turning from Topology to Dynamical Systems,  it was another Aha! moment.

Bill's main tool is of course the compactification of Teichm\"uller space by the projectified
space of measured foliations, yielding a space homeomorphic to a ball,
on which the mapping class group acts naturally. Therefore by Brouwer's fixed point theorem, each 
element has a fixed point in this compactification. The underlying geometrical nature of the fixed point
gives the classification.

The Orsay seminar on Thurston's work was very lively. The number of attendants was large.
Jean-Christophe Yoccoz who was just
starting graduate school told me that he attended it (I do not remember that, I really did not meet
him till the end of that academic year) and it left on him a lasting impression.

One of the main challenges during the lectures was the discussions with complex analysts
who had a compactification of Teichm\"uller space as a Euclidean ball by quadratic differentials.
The discussions were  driven by the belief that these two compactifications were the same.
It was a surprise when, sometime during the year, we learned that Steve Kerckhoff, then a
PhD student of Bill, showed that the two compactifications were distinct.
Of course, both compactifications are nowadays important, and they can be used to prove
the classification of elements in the mapping class group.

After that, the lamination point of view pervaded the subject. It was quite remarkable 
that Bill Thurston and Mike Handel were able
to show using laminations that the ideas of Nielsen that dated back to the 1930's potentially
contained the classification of elements of the mapping class group.
At that time, I was already getting back
to dynamics problems and  lost track of the subject.

Twenty years later laminations (not necessarily geodesic) came back to haunt me. 
There is hardly a day in my mathematical life without thinking about laminations.

In fact, about 1982, John Mather and Serge Aubry established the now so-called Aubry-Mather
theory for twist maps of the annulus. Although not usually expressed that way, the Aubry-Mather
set (or rather its suspension) is a lamination (not geodesic).
When Mather generalized these results to higher dimensions in the setting of Lagrangian systems,
the connection became much clearer. Aubry-Mather sets are foliated by 1-dimensional trajectories.
They are therefore laminations. Mather's graph theorem is in fact a proof that these laminations
are Lipschitz (the speed of the trajectory is a Lipschitz function of the point). This is a  crucial property
for geodesic laminations in dimension two, which follows in that case from a simple (hyperbolic)
planar argument.

In 1996, I discovered (like Weinan E and Craig Evans-Diogo Gomes) the relation between 
the Aubry-Mather theory and the viscosity solutions of the Hamilton-Jacobi equation.
The fact that I knew the lamination theory set up by Bill was definitely instrumental going deeper
in this relationship that keeps me still mathematically busy today.

I do not mention the work of Bill on Poincar\'e's conjecture and on holomorphic dynamics both of which had immense influence on several French mathematicians.
I personally had not been involved in these parts.

In the Fall of 1986, Universit\'e Paris-Sud
(Orsay) gave a Doctorat Honoris Causa to both Bill Thurston and Don Knuth.
At that time a plane ticket was a physical piece of paper that you actually needed to have 
to take the plane. Of course, the University President's staff bought an expensive ticket and they
were worried to send it by (regular mail): UPS, Fedex etc. with their overnight delivery were not
operating in France or at least not thought of. Anyway, I was planning to spend the 1986-87
academic year at IAS in Princeton, so one day Jean Cerf came to my office 
and asked me to deliver the ticket to Bill as soon as I would arrive, and to ask Bill to notify by fax the staff
of the University Presidency that the plane ticket has been delivered. It was very stressful:
I seem to remember that the price of the ticket was more than my monthly salary.
I arrived in Princeton late in the afternoon, hardly slept that night, first thing next day I ran
to Fine Hall, found Bill, delivered the ticket and followed him to the secretary's office to make sure
that the fax was sent. I felt much better afterwards.

I would like to end by mentioning New College in Sarasota (Florida) where Bill did his undergraduate studies
(John Smillie was also an undergraduate there). I think that the informal and congenial atmosphere at this wonderful
institution was instrumental in Bill's mathematical formation. The (apocryphal?) story I heard is that Bill spent four
years at New College essentially reading Fricke and Klein's book (in German!). I
 discovered New College during my years
at the University of Florida. After I returned to France to work at ENS Lyon, I visited several times New College trying (unsuccessfully) to attract some of the students to do their graduate studies in France. I was hoping that French mathematics
would return to the next Bill Thurston, at least a small part of what Bill gave us.

\section{Bill Abikoff}

 I spent the academic year 1976-77 as a Sloan Fellow at IH\'ES. I was hoping it would be a quiet place to work on Kleinian groups and indeed it was so quiet that I tiptoed in the halls so as not to disturb anyone. 
  
That all changed in the spring semester. Sullivan came back from the US and Gromov also arrived. Almost immediately, there were informal seminars in the hall with participants seated on the floor and sometimes shouting at each other. While in the States, Sullivan had proved that a compact metric manifold, of dimension unequal to four, admits a compatible Lipschitz structure. The theorem was unknown at IH\'ES except to the Director, Nico Kuiper, who had heard of it during a visit to the US. Kuiper asked Sullivan to lecture on that theorem. 
 
  A seminar at Orsay had already started trying to understand Thurston's work on 2-manifolds. It  was related to Teichm\"uller theory, and I was asked to lecture on Teichm\"uller's theorems. Bers arrived and he lectured on his complex-analytic proof of Thurston's result on the classification of mapping classes.
 
 I mentioned in the seminar the question of whether a change of basepoint in Bers' embedding of the Teichm\"uller space induces a map of Teichm\"uller space which extends to the Bers boundary. My interest in the problem was its consequence that the whole mapping class group extends to the boundary. Neither of these results is true and we later learned that Thurston used a geometry on the Teichm\"uller space which is quite different from that of Bers.
 
 By the spring, there were several seminars related to hyperbolic geometry at both IH\'ES and Orsay. Bill Harvey lectured on the \emph{curve complex} he had introduced; it is currently of great interest.
 
 Even by the end of the Orsay seminar, we didn't really understand Thurston's 2-manifold work. The issue was how to attach the boundary to Teichm\"uller space in a fashion that the mapping class group action extends.
 People like me were still thinking in terms of classical Teichm\"uller theory, and not in terms of hyperbolic geometry.   
  
Thurston had already moved on from surfaces to 3-manifolds. Siebenmann, who commuted several times that year between Orsay and Princeton came back with news about Thurston's \emph{bounded image theorem}, and told us that Thurston had announced a proof of the hyperbolization theorem for 3-manifolds. 

Marden had already shown that hyperbolic 3-manifolds, which arise in the context of classical Kleinian groups, are sufficiently large in the sense of Waldhausen; the hyperbolization theorem is a geometrization of the construction algorithm of Haken using, in the non-fibered case, the combination theorems due to Maskit.  People decided that we should forget about the planned lectures and concentrate on that. I outlined the proof for non-fibered Haken manifolds in a four hour marathon session.  The details of Thurston's ideas didn't even start to appear, in the notes prepared by Floyd and Kerckhoff, until two years later.

\section{David Fried}

I came to IH\'ES in the spring of 1977 to meet and work with Dennis Sullivan. I soon learned that many mathematicians in Paris, including Dennis, were obsessed with the new results of Thurston and that there was an active seminar devoted to his remarkable work on surfaces and 3-manifolds. This was a learning opportunity for me and I was pleased to play a small part in this seminar.

It began when Dennis spoke on a novel invariant, the Thurston norm $N$ of an oriented closed 3-manifold $M$. $N$ is a geometrically defined seminorm on the first cohomology of $M$ that takes integer values at each integral class $u$. Roughly speaking $N(u)$ is the minimal value of $\vert e(S)\vert$, where $S$ is  a closed aspherical surface in M associated to $u$ and $e(S)$ is its Euler characteristic. Dennis vigorously explained why $N$ was a seminorm. Using foliation theory, especially Thurston's thesis, he showed that $N(u) = \vert e(F)\vert$ when $u$ corresponds to a fibration of M over the circle with fiber $F$. Dennis described Thurston's theorem that a seminorm in finite dimensions with integer values at integer points must have a finite-sided unit ball but he did not remember the proof and the results were not yet in preprint form.

I found a proof, however, and soon found myself in Orsay presenting it to the seminar. This was the first talk I ever gave in France and I recall one incident from it fondly. Someone in the audience inquired whether a seminorm with rational values at rational points must also have a finite-sided unit ball. I admitted that I didn't know the answer and I returned to the blackboard. Another participant rose, pondered this delightful question, stroked his beard, wandered the room, chatted to himself, and began to use the far end of my blackboard for his scratchwork. No one seemed to find this odd and I happily carried on with my talk. 

I learned subsequently that the bearded thinker was Adrien Douady. At the next meeting of the seminar he presented his elegant counterexample: the norm $N$ on the Cartesian plane whose unit ball is the convex hull of the union of two unit discs with centers $(-1,0)$ and $(1,0)$. The unit sphere of $N$ meets each line through $(0,0)$ with rational slope in a rational point, so $N$ takes rational values at rational points.

I hope this suggests the fresh and open character of the Orsay seminar, which gave Thurston's work the close attention it deserved.
 
  \section{Dennis Sullivan}
 
\centerline{\it (A Decade of Thurston Stories)}

\medskip

\noindent {\bf First story.}\footnote{Editor's note: From an email Sullivan sent  to Athanase Papadopoulos, on April 27, 2019: ``I wrote these stories at one sitting soon after Bill passed away."} In December of 1971, a dynamics seminar  ended at Berkeley with
the solution to a thorny problem in the plane which had a
nice application in dynamics. The solution purported to move $N$
distinct points to a second set of  ``epsilon near" $N$  distinct points by
a motion which kept the points distinct and only moved while staying always
``epsilon prime near". The senior dynamicists in the front  row
were upbeat because the dynamics application up to then had only been
possible in dimensions at least three where this matching problem is
obvious by
general position. But now the dynamics theorem also worked in dimension two.

A heavily bearded long haired graduate student in the back of the room
stood up and said he
thought the algorithm of the proof didn't work. He went shyly to the
blackboard and drew two configurations of about seven points each and
started applying to these the method of the end of the lecture. Little
paths started emerging and getting in the way of other emerging paths
which to avoid collision had to get longer and longer. The algorithm
didn't work at all for this quite involved diagrammatic reason.
I had never seen such comprehension  and creative
construction of a counterexample done so quickly.
This combined with my awe at the sheer complexity of the geometry
that emerged.
\medskip

\begin{figure}[htbp]
\centering
\includegraphics[scale=.55]{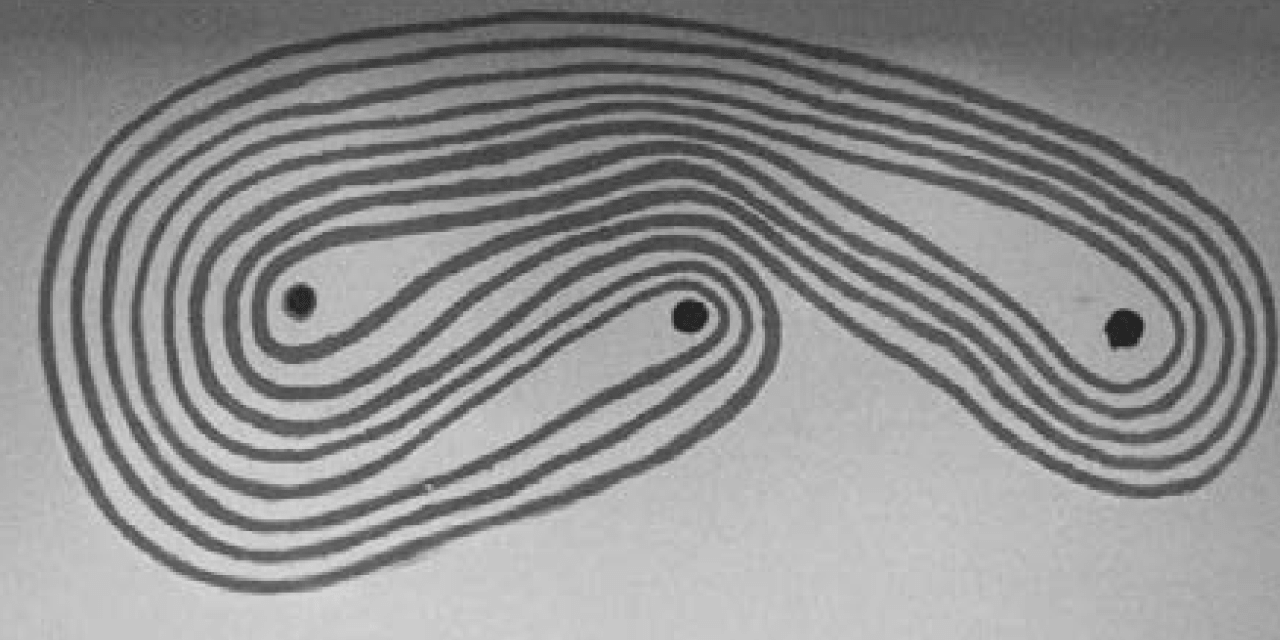}
\caption{The Berkeley wall curve painting by Thurston and Sullivan} \label{Berkeley}
\end{figure}

\noindent {\bf Second story.} A couple of days later the grad students invited me (I was
also heavily bearded  with long hair) to paint math frescoes on the
corridor wall  separating their offices from the elevator foyer.
While milling around before painting that same grad student came up to ask
``Do you think this is interesting to paint?" It was a complicated  smooth
one-dimensional object encircling three points in the plane. I asked 
``What is it?" and was astonished to hear ``It is a simple closed curve." I
said ``You bet  it's interesting!". So we proceeded to spend several
hours painting this curve on the wall. It was a great learning and
bonding experience.
For such a curve to look good it  has to be drawn in sections of short
parallel slightly curved strands (like the flow boxes of a foliation)
which are  subsequently smoothly spliced together.
When I asked how he got such curves, he said by successively applying to a
given simple curve a pair of Dehn twists along intersecting curves.
The ``wall curve painting", two meters high and four meters wide dated and signed,  lasted on that Berkeley
wall with periodic restoration for almost four decades before finally being
painted over a few years ago (see Figure \ref{Berkeley}).

\medskip

\noindent {\bf Third story.} That week in December 1971 I was visiting Berkeley from MIT
to give a
series of lectures on differential forms and the homotopy theory of
manifolds. Since foliations and differential forms were appearing
everywhere, I thought to use the one-forms that emerged in my story
describing the lower central series of the fundamental group to construct
foliations. Leaves of these foliations would cover graphs of maps of the
manifold to the nilmanifolds  associated to all the higher nilpotent
quotients of the fundamental group. These would generalize Abel's map to
the torus associated with the first homology torus.
 Being uninitiated in Lie theory I had asked all the differential geometers
at MIT and
Harvard about this possibility but couldn't make myself understood.
It was too vague/too algebraic.
I presented the discussion in my first lecture at Berkeley and to Bill
privately without much hope because of the  weird algebra/geometry mixture.
However the next day Bill came
with a complete solution and a full explanation. For him it was elementary
and
really only involved actually understanding the basic geometric
meaning of the Jacobi relation in the  Elie Cartan $d\circ d=0$ dual form.

In between the times of the first two stories above I  had spoken to my old
friend Mo Hirsch about  Bill Thurston who was working with Mo and was
finishing in his fifth year after an apparently slow start.  Mo or someone
else  told how  Bill's oral exam was a slight problem because when asked
for an example of the universal cover of a space Bill chose the surface of
genus two and started drawing awkward octogons with many [eight] coming
together at each vertex. This exposition quickly became an unconvincing
mess on the blackboard. I think Bill was the only one in the room who had
thought about this  nontrivial kind of universal cover.
 Mo then said ``Lately, Bill has started solving thesis level problems at
the rate of about one every month.'' Some years later I heard from Bill that
his first child Nathaniel didn't like to sleep at night so Bill was sleep
deprived ``walking the floor with Nathaniel'' for about a year of grad school.

 That week of math at Berkeley was life changing for me. I was very
grateful to be able to seriously appreciate the Mozart-like phenomenon I
had  been
observing; and I had a new friend.

Upon returning to MIT after the week in Berkeley I related my news to the
 colleagues there, but I think my enthusiasm was too intense to be
believed: ``I have just met the best graduate student I have ever seen or
ever expect to see." It was arranged for Bill to give a talk at MIT which
evolved into  a plan for him to  come to MIT after going first to IAS in
Princeton. It turned out he did come to MIT for just one year 1973-74.\footnote{That
year I visited IHES where I ultimately stayed for twenty odd years while Bill
was  invited back to Princeton, to the University.}

\medskip

\noindent {\bf Fourth story, IAS Princeton 1972--73.} When I visited the environs of Princeton from MIT in 72-73 I
had a chance to interact more with Bill. One day walking outside towards
lunch at IAS, I asked Bill what a horocyle was.
 He said ``you stay here" and he started walking away into the Institute
meadow.

After some distance he turned and stood still saying ``You are on the
circumference of a circle with me as center." Then he turned, walked much
further away, turned back and said something which I couldn't hear because
of the distance. After shouting back and forth to the amusement of the
members we realized he was saying the same thing ``You are on the
circumference of a circle with me as center". Then he walked even further
away, just a small figure in the distance and certainly out of hearing,
whereupon he turned and started shouting presumably the same thing again
and again. We got the idea what a horocycle was.

The next day, Atiyah asked some of us topologists if we knew if flat vector
bundles had a classifying space. (He had constructed some new characteristic
classes for such.) We knew it existed from Brown's theorem but didn't know
how to construct it explicitly. The next day, Atiyah said he asked Thurston
this question who did it by what was then a shocking construction: take the
Lie structure group of the vector bundle as an abstract group with the
discrete topology and form its classifying space.

Later, I heard about Thurston drawing Jack Milnor a picture proving any
dynamical pattern for any unimodal map appears in the quadratic family
$x\mapsto x^2 + c$. Since I was studying dynamics, I planned to spend a semester
with Bill at Princeton to learn about the celebrated Milnor--Thurston
universality paper that resulted from this drawing.

\medskip

\noindent {\bf Fifth story, Princeton University fall 1976.} I expected to learn about one-dimensional dynamics upon
arriving in Princeton in September 1976, but Thurston had already developed a
new theory of surface transformations. The first few days, he expounded on
this in a wonderful three hour  extemporaneous lecture at the Institute.
Luckily for me, the main theorem about limiting foliations was intuitively
clear because of the painstaking  Berkeley wall curve painting described
above.

At the end of that semester Bill told me he believed the mapping
torus of these carried hyperbolic metrics. When I asked why, he told me he
couldn't explain it to me because I didn't  understand enough differential
geometry. 

 A few weeks after, I left Princeton, with more time to work and without my
distractions. Bill essentially understood the proof of the
hyperbolic metric for appropriate Haken manifolds. The  mapping torus case
took two more years as discussed below.
During the semester  grad course  that Bill gave, the grad students
and I learned several key ideas:

 1) The quasi-analogue of ``hyperbolic geometry at infinity becomes
conformal geometry on the sphere at infinity"." (A notable memory here is the
feeling that Bill conveyed about really being inside hyperbolic space
rather than being outside and looking at a particular model. For me this
made a psychological difference.)

2)  We learned about the intrinsic geometry of convex surfaces outside the
extreme points: Bill came into class one day, and, for many minutes, he rolled a paper
contraption he had made around and around in the lecturer's table  without
saying a word until we felt the flatness.

 3) We learned  about the thick-thin decomposition of hyperbolic surfaces.
I remember how Bill drew a 50 meter long thin part winding all around the
blackboard near the common room, and suddenly everything was clear.
Including geometric convergence to the points of the celebrated  DM
compactification of the space of Riemann surfaces.

During that fall '76  semester stay at Princeton, Bill and I discussed
 understanding the Poincar\'e conjecture by trying to prove a general
theorem about
all closed three manifolds based on the idea that three is a relatively
small dimension.  We included in our little paper on ``canonical
coordinates" the sufficient for Poincar\'e Conjecture  
possibility that all closed
three manifolds carried conformally flat coordinate atlases.\footnote{This class is
closed under connected sum and contains many prime three manifolds.} However,
an undergrad, Bill Goldman, who was often around, disproved this a few years later for the
nilmanifold prime.\footnote{When looking for Mo Hirsch's current email, I noticed he had over over
200 descendants with a dozen coming from all but two students, Bill Goldman
with about 30 descendants and Bill Thurston with the rest.}   
We decided to try to spend a year together in the future.
\medskip

\noindent {\bf Meeting in the Alps,  spring 1978.}
 In the next period Bill developed limits of quasi-Fuchsian Kleinian groups
and pursued the mapping torus hyperbolic structure in Princeton while I
pursued the Ahlfors limit set measure problem in Paris. After about a year
Bill had made substantial positive progress (e.g., closing the cusp) and I had
made substantial negative progress (showing all known ergodic methods
coupled with all  known Kleinian group information were inadequate: there
was too much potential nonlinearity). We met in the Swiss Alps at the Plans-sur-Bex
 conference and compared notes. His mapping torus program was
positively finished but very
complicated while my negative information actually  revealed a rigidity
result extending Mostow's,
which allowed a simplification of 
Bill's proof.\footnote{See my Bourbaki report
on Bill's mapping torus theorem \cite{Sullivan}.}

\medskip

\noindent {\bf Sixth story, The Stony Brook meeting  summer 1978.} There was a big conference on Kleinian Groups at Stony Brook
and Bill was in attendance but not as a speaker. Gromov and I got him to
give a lengthy impromptu talk outside the schedule. It was a wonderful trip
out into the end of a hyperbolic 3-manifold, combined with convex
hulls, pleated surfaces and ending laminations \ldots
During the lecture Gromov leaned over and said watching Bill made him feel
like ``this field hadn't  officially started yet."

\medskip

\noindent {\bf Seventh story, Colorado June 1980 to August 1981.}  Bill and I shared the Stanislaw Ulam visiting chair at
Boulder and ran two seminars, a big one drawing together all the threads
for the full hyperbolic theorem and a smaller one on the dynamics of
Kleinian groups and dynamics in general.

All aspects of the hyperbolic proof passed in review with many grad
students in attendance.

 One day in the other seminar Bill was late
and Dan Rudolph was very energetically explaining in  just one hour  a new
shorter version of an extremely complicated  proof. The theorem promoted an
orbit equivalence to a conjugacy between two ergodic transformations if the
discrepancy of the orbit equivalence was controlled. The new proof was due
to a subset of the triumvirate Katznelson, Weiss and Ornstein and was
notable because it could be explained in one hour whereas the first proof
 took a mini-course to explain. Thurston  at last came in and  asked me to
bring him up to speed, which I did. The lecture continued to the end with
Bill wondering in loud whispers what the difficulty was and with me
shushing him out of respect for the context. Finally, at the end, Bill said
just imagine a bi-infinite string of beads on a wire with finitely many
missing spaces and just slide them all to the left say. Up to some standard
bookkeeping this gave a  new proof. Later that day an awestruck Dan
Rudolph said to me he never realized before then just how smart Bill
Thurston really was.

\medskip

\noindent {\bf Eighth story, La Jolla and Paris end of summer 1981.} 
The Colorado experience was very good, relaxing in the
Thurston seminar with geometry (one day we worked out the eight geometries
and another day we voted on terminology ``manifolded" or ``orbifold") and
writing  several papers of my own on Hausdorff dimension, dynamics and
measures on dynamical limit sets.

Later closer to Labor day I was flying from Paris to La Jolla to give a
series of AMS lectures on the  dynamics stuff when I changed the plan and
decided instead to try to expose the entire Hyperbolic Theorem ``for the
greater good" and as a self imposed Colorado final exam.  I managed to come
up with a one-page sketch while on the plane. There were to be two lectures
a day for four or five days. The first day  would be okay, I thought: just
survey things and then try to improvise for the rest, but I needed a stroke
of luck. It came big time.

 There is a nine hour time difference between California and Paris and the
first day I awoke  around midnight  local time and went to my assigned
office to prepare. After  a few hours I had generated many questions and
fewer answers about the  hyperbolic argument.
I noticed a phone on the desk that miraculously allowed long distance
calls and by then it was around 4:00 a.m. California time and  7:00 a.m. in
Princeton.  
I called Bill's house, and he answered. I posed my questions. He gave quick
responses, I took notes, and he said call back after he dropped kids at
school and got to his office. I gave my objections to his answers around
 9:30 a.m. his time and he  responded more fully. We ended up with various
alternate routes that all in all covered  every point. By 8:00 a.m. my time I
had a pair of lectures prepared.  The  first day went well: lecture/lunch/beach/swim, second lecture, dinner then goodby to colleagues and
back to bed. This took some discipline but as viewers of the videos will see
the audience was formidable (Ahlfors, Bott, Chern, Kirby,
Siebenmann, Edwards, Rosenberg, Freedman, Yau, Maskit, Kra, Keen, Dodziuk,\ldots)
and I was motivated.

 Bill and I  repeated this each day, perfecting the back and forth so that
by 8:00 a.m. California time
each day, I had  my two lectures prepared and they were  getting the job
done. The climax came when presenting  Bill's delicious argument that
controls the length of a geodesic representing the branching locus of a
branched pleated surface by the dynamical rate of chaos or entropy  created
by the geodesic flow on the intrinsic surface.  One knows that this is
controlled by the area growth of the universal cover of the branched
surface which by negative curvature is controlled by the volume growth  of
the containing hyperbolic three space QED. There was in addition Bill's
beautiful example showing the estimate was qualitatively sharp. This
splendid level of lecturing was too much for Harold Rosenberg, my astute
friend from Paris, who was in the audience. He came to me afterwards and
asked frustratingly ``Dennis, do you  keep Thurston locked up in your
office upstairs?'' The lectures were taped by Michael Freedman and I have
kept my lips sealed until now. The taped (Thurston)-Sullivan lectures are available online.\\

\noindent {\bf Ninth story, Thurston in Paris fall 1981.} Bill visited me in Paris and I bought a comfy sofa bed for
my  home office where he could sleep.
He politely asked what would I have talked about had I not changed plans
for the AMS lectures, and in particular what had I been doing in detail in
Colorado beyond the hyperbolic seminar.
There were about six papers to tell him about. One of the most appealing
ideas I had learned from him.
Namely the visual  Hausdorff $f$-dimensional
measure of an appropriate set on
the sphere at infinity, as viewed from a point inside, defines a positive
eigenfunction for the hyperbolic 3-space Laplacian with eigenvalue
$f(2-f)$.

I started going through the ideas and statements. I made a statement and he
either  immediately gave the proof or I gave the idea of my proof. We went
through all the theorems in the  six papers in one session with either him
or me giving the proof.
There was one missing result that the bottom eigenfunction when $f$ was $> 1$
would be represented by a normalized eigenfunction whose square integral
norm was estimated by the volume of the convex core. Bill lay back for a
moment on his sofa bed, his eyes closed, and immediately proved the missing
theorem. He produced the estimate by diffusing geodesics transversally and
averaging.

Then we went out to walk through Paris from Porte d'Orl\'eans to Porte de
Clignancourt. Of course we spoke so much about mathematics that Paris was essentially
forgotten, except maybe the simultaneous view of Notre-Dame and the
Conciergerie as we crossed over the Seine.\\

\noindent {\bf Tenth story, Princeton-Manhattan 1982-83.} I began splitting time between IHES and the CUNY grad center
where I started a  thirteen year long Einstein chair seminar on dynamics
and quasiconformal homeomorphisms (which changed then to quantum objects in
topology) while Bill continued developing a cadre of young geometers to
spread the beautiful ideas of  negatively curved space.
Bill delayed writing a definitive text on the hyperbolic proof in lieu of
letting things develop along many opening avenues\footnote{I watched recently with great pleasure the unfolding of the ingenious
proof by Kahn and Markovic of the  hyperbolic subsurface conjecture. As
each step was revealed I remembered when some key/crucial aspect of an
analogous device was introduced by Bill more than thirty years before and
then later taught to his proteges at Princeton.} by his increasingly
informed cadre of younger/older geometers.
He wanted to avoid what happened when his basic papers on foliations
``tsunamied" the field in the early 1970s.

Once we planned to meet in Manhattan  
to discuss
holomorphic dynamics in one variable and its analogies with hyperbolic
geometry and Kleinian groups  that I had been preoccupied with.
We were not disciplined and began talking about other things at the
apartment, and finally got around to our agenda about thirty minutes before
he had to leave for his train back to Princeton. I sketched the general
analogy: Poincar\'e limit set, domain of discontinuity, deformations,
rigidity, classification, Ahlfors finiteness theorem, the work of
Ahlfors--Bers,\ldots to be compared with
Julia set, Fatou set, deformations, rigidity, classification, non-wandering
domain theorem, the work of Douady--Hubbard,\ldots
which he perfectly and quickly absorbed until he had to leave for the train.
Two weeks later we heard about his reformulation of a holomorphic dynamical
system as a fixed point on Teichm\"uller space analogous to part of his
hyperbolic theorem. There were many new results  including those of Curt
McMullen some years later and
 the subject of holomorphic dynamics was raised to another higher level.

    \medskip
    
\noindent {\bf Postscript.} Thurston and I met again at Milnor's 80th fest at Banff after  essentially
thirty years and picked up where we had left off. I admired his checked
green shirt the second time it appeared and he  presented  it to me the
next day. We  promised to try to attack together a remaining big hole in
the Kleinian group/holomorphic dynamics dictionary: ``the invariant line
field conjecture". It was a good idea but unfortunately turned out to be
impossible.
At the same conference, I recall a comment  whispered by Bill  who sat next to me during a talk by Jeremy Kahn about the Kahn--Markovich   proof  of the  Subsurface Conjecture from decades before.   Bill  whispered : ``I missed the `offset' step."\footnote{In the Kahn--Markovich  proof, one glues up all possible ideal triangles building many immersed surfaces, and uses ergodic theory of the natural actions on this space.
The offset idea is to glue the triangles not at their midpoints but after offsetting by a fixed  large amount.
This  prevents  missing everything  in the limit. This step was the one that was missed before by Bill.}

\bigskip

\noindent {\it Acknowledgements.} We would like to thank  Annette A'Campo Neuen, Vincent Colin, Yi Huang, Silvio Levy and Karen Vogtmann for their helpful comments and corrections on these articles, Harold Rosenberg for sending us a copy of the letter from Milnor to Thurston that is included here and Fran\c cois Tisseyre who kindly provided photos from the Clay conference in Paris.

\bigskip

\noindent  Fran\c cois Laudenbach, Laboratoire Jean Leray, Universit\'e de Nantes,  2  rue de la Houssini\`ere, F-44322 Nantes Cedex 3, francois.laudenbach@univ-nantes.fr; 
 
 \medskip
 
\noindent Athanase Papadopoulos,  Chebyshev Laboratory at St. Petersburg State University, Russia, 199178, St. Petersburg,
14th line of the Vasilyevsky Island, house 29B, and Institut de Recherche Mathématique Avancée,
7 rue René Descartes
67084 Strasbourg Cedex France, 
papadop@math.unistra.fr

\end{document}